\documentclass[11pt]{article}
\usepackage{amsfonts, amssymb, eucal, amsmath, eufrak}
\usepackage{graphicx}
\usepackage{srcltx}
\usepackage{amsmath}
\usepackage{amsthm}
\usepackage{amssymb}
\usepackage{makeidx}
\usepackage{subfloat}
\usepackage{subcaption}
\usepackage{algorithmic}
\usepackage{xcolor}
\usepackage{enumerate}
\usepackage{thmtools}
\usepackage{bm}
\usepackage{url}

\newtheorem{alg}{Algorithm}[section]
\newtheorem{lemma}{Lemma}[section]

\newtheorem{theorem}{Theorem}[section]

\newtheorem{remark}{Remark}[section]

\declaretheoremstyle[
bodyfont=\normalfont\itshape,
headformat=\NAME\NUMBER  
]{}
\declaretheorem[style=nospacetheorem,name=Assumption ]{manualtheoreminner}
\newenvironment{assumption}[1]{%
	\renewcommand\themanualtheoreminner{#1}%
	\manualtheoreminner
}{\endmanualtheoreminner}

\newcommand{\R}{\mathbb R}
\newcommand{\Rn}{\mathbb R^n}

\newcommand{\Rm}{\mathbb R^m}

\newcommand{\N}{\mathbb N}

\newcommand{\tr }{^\top}
\DeclareMathOperator*{\eig}{eig}

\DeclareMathOperator*{\dia}{diag}
\DeclareMathOperator*{\atan}{atan2}

\DeclareMathOperator*{\dist}{dist}

\newcommand{\calN}{\mathcal N}

\newcommand{\calK}{\mathcal K}

\newcommand{\xbf}{\mathbf x}
\newcommand{\rbf}{\mathbf r}
\newcommand{\zbf}{\mathbf z}

\newcommand{\gbf}{\mathbf g}

\newcommand{\dbf}{\mathbf d}

\newcommand{\ybf}{\mathbf y}

\newcommand{\Rbf}{\mathbf R}

\newcommand{\mmu}{m_{\mu}}
\newcommand{\Mmu}{M_{\mu}}

\newcommand{\lr}[1]{\left(#1\right)}

\newcommand{\REV}[1]{#1}

\newcommand{\be}{\begin{equation}}
	\newcommand{\ee}{\end{equation}}
\def\qed{$\square$}

\begin {document}
\title{Parallel Inexact Levenberg-Marquardt Method for Nearly-Separable Nonlinear Least Squares}
\author{Lidija Fodor\footnotemark[1]
	Du\v{s}an Jakoveti\'c\footnotemark[1],  Nata\v{s}a Kreji\'c\footnotemark[1],   Greta Malaspina\footnotemark[2]
}
\date{}
\maketitle
\footnotetext[1]{Department of Mathematics and Informatics, Faculty of Sciences, University of Novi Sad, Trg DositejaObradovi\'ca 4, 21000 Novi Sad, Serbia, Email: lidija.fodor@dmi.uns.ac.rs, dusan.jakovetic@dmi.uns.ac.rs, natasak@uns.ac.rs.  The work of Jakoveti\'c and Kreji\'c  was supported by the Science Fund of the Republic of Serbia, GRANT No 7359, Project title - LASCADO.}
\footnotetext[2]{Dipartimento di Ingegneria Industriale, Università degli studi di Firenze, Viale G.B. Morgagni 40,50134 Firenze, Italy, Member of the INdAM Research Group GNCS, Email: greta.malaspina@unifi.it. The work of Malaspina is partially supported by Partenariato esteso FAIR ``Future Artificial Intelligence Research'' SPOKE 1 Human-Centered AI. Obiettivo 4, Project ``Mathematical and Physical approaches to innovative Machine Learning technologies (MaPLe)''.}
\begin{abstract}
	Motivated by localization problems such as cadastral maps refinements, we consider a generic Nonlinear Least Squares (NLS) problem of minimizing an aggregate squared fit across all nonlinear equations (measurements) with respect to the set of unknowns, e.g., coordinates of the unknown points' locations. In a number of  scenarios, NLS problems exhibit a nearly-separable structure: the set of measurements can be partitioned into disjoint groups (blocks), such that the unknowns that correspond to different blocks 
	are only loosely coupled.  
	We propose an efficient parallel method, termed Parallel Inexact Levenberg Marquardt (PILM), to solve such generic large scale NLS problems. PILM builds upon the classical Levenberg-Marquard (LM) method, with a main novelty in that the nearly-block separable structure is leveraged in order to obtain a scalable parallel method. Therein, the problem-wide system of linear equations that needs to be solved at every LM iteration is tackled iteratively. At each (inner) iteration, the block-wise systems of linear equations are solved in parallel, while the problem-wide system is then handled via sparse, inexpensive inter-block communication. We establish strong  convergence guarantees of PILM that are analogous to those of the classical LM; provide PILM implementation in a master-worker parallel compute environment; and demonstrate its efficiency on huge scale cadastral map refinement problems.

	\noindent{Keywords}: Distributed optimization; sparse nonlinear least squares; Inexact Levenberg-Marquardt method; nearly separable problems; localization problems; cadastral maps refinement.\end{abstract}

\section{Introduction}

Localization is a fundamental problem in many applications such as Network Adjustment \cite{netadj}, Bundle Adjustment \cite{bundleadj}, Wireless Sensor Network Localization \cite{WSN}, target localization \cite{Claudia2023}, source localization \cite{AmirBeck}, etc. In localization-type problems, the goal is typically to determine unknown locations (e.g., positions in a 2D or a 3D space) of a group of points (e.g., entities, sensors, vehicles, landmarks) based on (often highly nonlinear and noisy) measurements that involve unknown positions corresponding to only a few entities; these measurements may be, e.g., point to point distance measurements, angular measurements involving three points, etc. Mathematically, a broad class of  localization-type applications may be formulated as \emph{Nonlinear Least Squares} (NLS) problems (see, e.g., \cite{JOGO3,JOGO4,JogoLOVO}). Specifically, we are interested in the following formulation:
\begin{equation}\label{LS}
	\min_{\xbf\in\R^n} F(\xbf),\enspace F  (\xbf) = \frac{1}{2}\sum_{j=1}^m r_j(\xbf)^2=\frac{1}{2}\| \Rbf(\xbf)\|^2.  
\end{equation}
Here, for every $j=1,\ldots,m$, $r_j: \R^n \to \R $, $ \Rbf(\xbf) = \left(r_1(\xbf),\dots,r_m(\xbf)\right)\tr \in\Rm$ is the vector of residuals, where $r_j$ stands for the residual function that corresponds to the $j$-th measurement; and $ F  $ is the aggregated residual function. We refer to Figure \ref{fig:measurements} and equations \eqref{eqn-measurement-1}--\eqref{eqn-measurement-3} ahead to illustrate the quantities $r_j$ (x) for point-to-point distance, angle, and point-line distance measurements that arise with localization-type problems. We assume that the problem is large scale, in the sense that both the dimension $n$ (total number of unknown per-coordinate locations) and the number of functions (measurements) $m$ are large. Furthermore, we assume that the problem is nearly separable, as detailed  below. 
While we are primarily motivated by localization-type problems, 
the method proposed here, as explained ahead, applies to generic NLS problems.

There has been a significant body of work that considers various localization-type problems of the form \eqref{LS} or of some special case of \eqref{LS} \cite{bundleadj,Claudia2020,netadj}. Unlike our formulation that keeps a general residual structure, many works focus on a specific type of residual functions, e.g., point to point distance only \cite{Claudia2015}, or point to point distance and angular measurements \cite{Claudia2023}. Typically, as it is the case here, the corresponding formulations \eqref{LS} are non-convex. Various approaches have been proposed to solve these formulations, including convex approximations of the objective \cite{Claudia2015}, linearization of measurements like in extended Kalman filter approaches (EKF) \cite{EKF}, and a direct handling of the non-convex problem, \cite{DistributedLMeurasip, Ers, SLM}, as it is done here. The corresponding optimization solvers may be centralized, e.g.,  \cite{AmirBeck}, or distributed, e.g., \cite{Claudia2023}.

Unlike most of existing studies, our main motivation is to efficiently solve large-scale problems \eqref{LS} that exhibit a nearly-separable structure. The concept of near-separable is introduced in \cite{SLM}. It roughly means that the set of measurements may be partitioned into disjoint groups, such that the variables involved in each of the measurement groups are only loosely coupled. The nearly separable structure of \eqref{LS} arises in many applications. One specific scenario that motivates the nearly separable structure is optimization-based improvement of cadastral maps \cite{franken, heuvel}. That is, given an existing (low-accuracy) cadastral map and a set of accurate field surveyor measurements, the goal is to improve the accuracy of the map by making it consistent with the available measurements. Since the map typically represents a whole country or large region, while surveyor measurements only involve points that are relatively close to each other (e.g. in the same neighborhood), the resulting problem is often nearly separable.

The paper's contributions are as follows. We propose an efficient parallel method to solve generic large scale problems \eqref{LS} with a nearly-separable structure, termed PILM (Parallel Inexact Levenberg Marquardt method). While this paper is mainly motivated by localization-type problems, PILM is suitable for arbitrary NLS problems. In the context of localization problems, PILM handles arbitrary types of measurements such as, e.g., point to point distances, angular measurements, bundle measurements, etc. The proposed approach builds upon the classical Levenberg-Marquard (LM) method for solving NLS. However, a main novelty here is that the nearly block separable structure is leveraged in order to obtain a scalable parallel method. In more detail, the problem-wide system of linear equations that needs to be solved at every LM iteration is tackled iteratively, where block-wise systems of linear equations are solved in parallel, while the problem-wide system of equations is then handled via sparse communication across the multiple subsystems. 

Specifically, PILM's implementation assumes a master-worker computational infrastructure, where the master orchestrates the work, while each worker possesses data that corresponds to one of the sub-problems. Then, the costly solving of the system of linear equations at each iteration is delegated to the worker nodes working in parallel, while the master node aggregates local solutions from the workers and  performs the global update. More precisely, at each iteration, a search direction is computed by approximately solving the linear system that arises at each iteration of LM in a distributed manner, using a fixed point strategy that relies on the partition of the variables induced by the near-separable problem structure.

We prove that PILM, combined with a nonmonotone line search strategy, achieves global convergence to a stationary point of \eqref{LS}, while for the full step size, we prove local convergence, with the convergence order depending on the choice of the parameters of the method. We see that, aside from near-separability, the required assumptions are the same as for the convergence analysis of the classical LM method. In other words, we prove convergence of the novel parallel method under a standard set of LM assumptions. \\
{ The nonlinear least squares problem \eqref{LS} that we consider is in
	general non-convex and very challenging to solve globally. However,
	interestingly, for certain special cases of \eqref{LS}, despite the problem’s
	nonconvexity, a global solution to the problem can be found. 
	In \cite{Beck} the authors considers \eqref{LS} when only range-based
	measurements are included and show that, under some
	mild conditions, the standard semidefinite programming (rank-1)
	relaxation of the original problem is exact, which enables them to achieve global optimality. Reference \cite{Rabbat}
	considers a further simplified version of \eqref{LS} with multiple range-based
	measurements and a single unknown target location. For this formulation,
	the authors study a normalized incremental gradient method and show that, despite the problem’s non-convexity, the method converges to the global optimum under some mild conditions. In this work, our goal is to devise an efficient parallel algorithm for a very general formulation of \eqref{LS}, hence guaranteeing global optima is extremely difficult. We therefore prove analytically global
	convergence to a stationary point. Moreover, we provide several
	numerical studies that demonstrate that the stationary point to which
	the algorithm converges corresponds to accurate
	localization solutions for the studied problems.}\\

Numerical examples on large scale cadastral map problems (with both the number of unknowns and the number of measurements on the order of a million) 
demonstrate efficiency of PILM and show that it compares favorably with respect to the classical LM method for large scale and nearly separable settings. An open source parallel implementation of PILM  is available at {  \url{https://github.com/lidijaf/PILM}.}

From the technical standpoint, many modifications of the classical LM scheme have been proposed in literature. Typically, the goal is to retain convergence of the classical LM while relaxing the assumptions on the objective function,  and to improve the performance of the method. 
In \cite{loc1, loc2, loc3} the damping parameter is defined as a multiple of the objective function. With this choice of the parameter local superlinear or quadratic convergence is proved under a local error bound assumption for zero residual problems, while global convergence is achieved by employing a line search strategy. In \cite{santos:global} the authors propose an updating strategy for the parameter that, in combination with Armijo line search, ensures global convergence and $q$-quadratic local convergence under the same assumption of the previous papers. In \cite{santos:local} the non-zero residual case is considered and a LM scheme is proposed that achieves local convergence with order depending on the rank of the Jacobian matrix and of a combined measure of nonlinearity and residual size around stationary points.
In \cite{inexact1} an Inexact LM is considered and local convergence is proved under a local error bound condition. In \cite{Bellavia} the authors propose an approximated LM method, suitable for large-scale problems, that relies on inaccurate function values and derivatives. A sequential method for the solution of the sparse and nearly-separable problems of large dimension  problems was proposed in \cite{SLM}. 
Unlike existing works, we develop a novel parallel method adapted to nearly separable structures, ensuring its theoretical guarantees and providing its parallel implementation with a validated efficiency on large-scale problems.

It is also worth noting that LM has been applied for solving distributed sensor network localization problems in~\cite{DistributedLMeurasip}. Therein, the authors are concerned with point to point distance measurements only. The method proposed therein is different than ours: it uses message passing and utilizes communication across a tree of cliqued sensor nodes. 
In contrast, we use an iterative fixed point strategy to leverage the block separable structure assumed here, and we also provide a detailed convergence analysis of the proposed method. 

This paper is organized as follows.
In Section \ref{sec:PLS_problem} we formalize the near-separability assumption, describe the induced block-partition of the LM equation and present PILM. In Section \ref{sec:PLS_convergence} we carry out the convergence analysis, while in Section \ref{sec:PLS_results} we present a set of numerical results to investigate the performance of the method. 

The notation we use in the paper is as follows. Vectors are denoted by boldface letters, $ \|\cdot\| $ is the 2-norm for both vectors and matrices, $\mathbb N_0:=\mathbb N\cup\{0\}$, where $\mathbb N$ denotes the set of natural numbers. Given a matrix $A\in\R^{n\times n}$ we denote with $\lambda_{\min}(A)$ and $\lambda_{\max}(A)$ the smallest and largest eigenvalue of $A$, respectively, in absolute value. That is, we take
$\lambda_{\min}(A) = \min\{|\lambda_i|\ :\ \lambda_i\ \text{eigenvalue of}\ A\}$,
and $\lambda_{\max}(A)$ is defined analogously.

\section{Distributed Inexact LM method }\label{sec:PLS_problem}
The problem we consider is stated in (\ref{LS}). A standard iteration of LM method is, for a given iteration $ \xbf^k $
\begin{equation*}
	\xbf^{k+1} = \xbf^k+ \dbf^k,
\end{equation*}
where $\dbf^k\in\R^n$ is the solution of
\begin{equation}\label{LM_system}
	\left( (J^k)\tr  J^k+\mu_k  I \right)\dbf^k = - (J^k)\tr  \Rbf ^k,
\end{equation}
where $ J^k =  J(\xbf^k)\in\R^{m\times n}$ denotes the Jacobian matrix of $ \Rbf ^k =  R (\xbf^k)$ and $\mu_k>0$ is a scalar. 
When $n$ is very large solving \eqref{LM_system} at each iteration of the method may be prohibitively expensive. In the following we propose an Inexact Levenberg-Marquardt method that relies on the near-separability of the problem to define a fixed-point iteration for the solution of the linear system \eqref{LM_system}. Such a method is suitable for the server/worker framework as it relies on the fixed point iterations that can be efficiently carried out in parallel, with modest communication traffic,  due to the near-separable structure of the problem. Let us first recall the near-separable property as introduced in \cite{SLM}. The description bellow, although a bit lengthy is included for the sake of completeness.
Let us define $\mathcal{I} = \{1,\dots,n\}$ and $\mathcal{C} = \{1,\dots,m\}$. Given a partition $I_1,\dots,I_{K}$ of $\mathcal I$ we define the corresponding partition of $\mathcal{C}$ into $E_1,\dots, E_{K}$ as follows:
\begin{equation}\label{def_Es}
	\begin{aligned}
		&E_s = \{j\in\mathcal{C} | r_j \text{ only depends on variables in } I_s\},\ s=1,\ldots,K\\
		&\widehat E = \mathcal{C} \setminus\bigcup_{i=s}^{ K }E_i.
	\end{aligned}
\end{equation}
In other words each of the subsets $E_s$ contains the indices corresponding to residual functions that only involve variables in $I_s$. The    indices of residuals that involve variables belonging to different subsets $I_s$ are gathered in $\widehat E$. If there exist $ K \geq2$ and a partition $\{I_s\}_{s=1,\ldots,K }$ of $\mathcal I$ such that $\widehat E = \emptyset$ we say that the problem \eqref{LS} is separable,  while we say that it is \emph{nearly-separable} if there exist $ K \geq2$ and a partition $\{I_s\}_{s=1,\ldots,K }$ of $\mathcal I$ such that the cardinality of $\widehat E$ is small with respect to the cardinality of $\bigcup_{s=1}^{ K }E_s.$ Let us explain with more details what happens with \eqref{LM_system} in the case of separable and near-separable problems.

Assuming that the partitions of the variables and residuals,  $\{I_s\}_{s=1}^K$ and ${ \{E_s\}_{s=1}^K} ,\ \widehat E,$ 
we define $\xbf_{s}\in\R^{n_s}$ as the vector of the variables in $I_s. $ Here  $n_s = |I_s| $ denotes the cardinality of $I_s$. Let us introduce the following functions
\begin{equation}\label{local_functions}
	\begin{aligned}
		&\Rbf_s(\xbf_s) := { (r_j(\xbf_s))_{j\in E_s}}, \phantom{space} &\rho(\xbf) := (r_j(\xbf))_{j\in \widehat E}\\
		&F_s(\xbf_s) :=\|\Rbf_s(\xbf_s)\|^2 \phantom{space} &\Phi(\xbf) := \|\rho(\xbf)\|^2
	\end{aligned}
\end{equation}
so that for every $s=1,\ldots,K $, $\Rbf_s:\ \R^{n_s} \to \R^{|E_s|} $ is the vector of residuals involving only variables in $I_s$, while $\rho: \R^n \to \R^{|\widehat E|} $ is the vector of residuals in $\widehat E$ and $F_s$, $\Phi$ are the corresponding local aggregated residual functions.  Clearly $ \sum_{s=1}^K n_s = n $ as  $\{I_s\}_{s=1}^K$ is a partition of $\mathcal{I}. $  
Now, problem \eqref{LS} is equivalent  written to
\begin{equation}\label{LS_sep}
	\begin{aligned}
		&\min_{\xbf\in\R^n}\left( \Phi(\xbf)+ \sum_{s=1}^{K } F_s(\xbf_s)\right) = \min_{\xbf\in\R^n}\left( \|\rho(\xbf)\|^2 + \sum_{s=1}^{K } \|\Rbf_s(\xbf_s)\|^2 \right)
	\end{aligned}
\end{equation}
The above expression implies that if the set $\widehat E$ is empty, i.e., if the problem is separable, then $\Phi(\xbf) = 0$ for every $\xbf\in\Rn. $ In this case one can solve  \eqref{LS_sep}  by solving $ K $  independent problems defined by 
\begin{equation}\label{LS_subp}
	\begin{aligned}
		&\min_{\xbf_s\in\R^{n_s}}F_s(\xbf_s) = \min_{\xbf_s\in\R^{n_s}} \|\Rbf_s(\xbf_s)\|^2\ \text{ for}\ s=1,\ldots, K .
	\end{aligned}
\end{equation}
If the set $\widehat E$ is not empty then the problem is not separable, implying that in general $\Phi$ is not equal to zero.  

Given the partitions $\{I_s\}_{s=1,\ldots,K }$ and $\{E_s\}_{s=1,\ldots,K }$ of $\mathcal I$ and  $\mathcal{C}, $ for simplicity we assume that the vectors $\xbf$ and $\Rbf$ are reordered according to these partitions.  In other words,   for $ \xbf \in \R^n $
$$
\xbf = \left(\begin{matrix}
	\xbf_1\tr, \ldots, 
	\xbf_{ K }{\tr}
\end{matrix}\right){\tr}, 
\phantom{spa}
\Rbf (\xbf) = \left(\begin{matrix}
	\Rbf_1(\xbf_1){\tr}, \ldots, 
	\Rbf_{ K }(\xbf_{ K }){\tr}, 
	\rho(\xbf){\tr}
\end{matrix}\right){\tr}
$$
Let us now look at the structure of the Jacobian. Denoting by $J_{s\Rbf_j}$  and  $ J_{s \rho} $ the partial derivatives,
$$ J_{s\Rbf_j} = \frac{\partial \Rbf_j}{\partial \xbf_s}, \; J_{s \rho} = \frac{\partial \rho}{\partial \xbf_s}, s=1,\ldots,K,  
$$
the Jacobian is a block matrix
\begin{equation*}
	J(\xbf) = \left(\begin{matrix}
		J_{1\Rbf_1}(\xbf_1)&       &   &   0    \\
		&J_{2\Rbf_2}(\xbf_2)&   &       \\
		&			  & \ddots&       \\
		0   &       &   &J_{ K \Rbf_{ K }}(\xbf_{ K })\\
		J_{1\rho}(\xbf)& J_{2\rho}(\xbf) & \dots &J_{ K \rho}(\xbf)
	\end{matrix}\right).
\end{equation*}

The diagonal blocks $J_{i\Rbf_i}(\xbf_i), i=1,\ldots, K$ depend only on variables in the set $ I_i $ and involve the residual functions in $ E_i.$ Therefore, assuming that we have $K$ worker computational nodes, such that node $i$ holds the portion of the dataset relative to the functions in $E_i$, we have that each of the nodes can compute one of the diagonal blocks of the Jacobian. Communication is only needed for the last row of blocks.
From this structure of $ \Rbf$ and $ J$ we get the corresponding block structure of the gradient $ \gbf(\xbf) =  J(\xbf)\tr  \Rbf (\xbf)$ and the matrix $ J(\xbf)\tr  J(\xbf)$:
\begin{equation} \label{7} 
	\gbf(\xbf)\tr= \left(\gbf_1\tr(\xbf), \gbf_2\tr(\xbf), \ldots, \gbf_{ K }\tr(\xbf)\right),
\end{equation}
\begin{equation}\label{JHB}
	J(\xbf)^\top  J(\xbf) = \left(\begin{matrix}
		P_1(\xbf)& B_{12}(\xbf)&\dots   &   B_{1 K }(\xbf)    \\
		B_{21}(\xbf)  &P_2(\xbf)& \ddots   &  \vdots     \\
		\vdots  &\ddots  & \ddots&  B_{ K -1 K }  (\xbf)    \\
		B_{ K 1}(\xbf)  &\dots  & B_{ K  K -1}(\xbf)   &P_{ K }(\xbf)\\
	\end{matrix}\right),
\end{equation}
where, for $s,i,j=1,\dots, K$
\begin{equation}\label{HandB_1}
	\begin{aligned}
		&\gbf_s(\xbf)= J_{s\Rbf_s}(\xbf_s)\tr \Rbf_s(\xbf_s) + J_{s\rho_s}(\xbf)\tr \rho(\xbf) \in\R^{n_s}\\
		&P_{s}(\xbf) = J_{s\Rbf_s}(\xbf_s)\tr J_{s\Rbf_s}(\xbf_s) + J_{s\rho}(\xbf)\tr J_{s\rho}(\xbf)\in\R^{n_s\times n_s}\\
		&B_{ij}(\xbf) = J_{i\rho}(\xbf)\tr J_{j\rho}(\xbf)\in\R^{n_i\times n_j} .
	\end{aligned}
\end{equation}
With this notation we can write
$$J(\xbf)^\top  J(\xbf) = P(x)+B(x),$$
where $ P(x) $ is the block diagonal matrix with diagonal blocks given by $P_s(\xbf)$ for $s=1,\dots, K $ and $ B(\xbf)$ is the block partitioned matrix with diagonal blocks equal to zero and off-diagonal blocks equal to $B_{ij}(\xbf)$. 

The algorithm we introduce here is motivated by the special structure of the Jacobian called near-separability property. We are now in position to define it formally. 
\begin{assumption}{E1} \label{ass:LS} There exists a constant $ C_B > 0 $ such that for all $ \xbf \in \R^n $
	\begin{equation} \label{M}
		\| B(\xbf)\| \leq C_B \| J(\xbf)\tr  J(\xbf)\|. 
	\end{equation}
\end{assumption} 
We notice that $B(\xbf) $ is a submatrix of $  J(\xbf)\tr  J(\xbf)$ and there is no upper bound over the magnitude of $C_B$, so the assumption above is not restrictive.

The linear system \eqref{LM_system} at iteration $k$ can be rewritten as
\begin{equation} \label{systemFP}
	( P^k +\mu_k  I+ B^k)\dbf^k = -\gbf^k
\end{equation}
where $\gbf^k = \gbf(\xbf^k),\ P^k = P(\xbf^k)$ and $B^k = B(\xbf^k)$ are defined above
with $\xbf = \xbf^k$. 

Consider the sequence $\{\ybf^l\}$ generated as follows
\begin{equation}\label{fixedpoint}
	\begin{cases}
		\ybf^1 = -( P^k+\mu_k I)^{-1}\mathbf g^k & \\
		\ybf^{l+1} = -( P^k+\mu_k I)^{-1}\left(\mathbf g^k+ B^k\ybf^l\right) & l\geq 1.
	\end{cases}
\end{equation}
Notice that $P(x)$ is a block diagonal matrix with $s$-th diagonal block given by $$P_{s}(\xbf) = J_{s\Rbf_s}(\xbf_s)\tr J_{s\Rbf_s}(\xbf_s) + J_{s\rho}(\xbf)\tr J_{s\rho}(\xbf)\in\R^{n_s\times n_s}.$$ This implies that $P(\xbf)$ is positive semi-definite for every $\xbf\in\Rn$, and therefore the matrix $( P^k+\mu_k I)$ is invertible for any $\mu>0.$
The equations above define a fixed-point method for { computing} the solution of \eqref{systemFP}. From the theory of fixed point methods we know that if $$\| ( P^k+\mu_kI)^{-1}B^k\|<1$$
then the sequence $\{\ybf^l\}$ converges to the solution $\dbf^k$ of \eqref{systemFP}.\\
Moreover, denoting with $\rbf^l_k$ the residual in the linear system at the $l$-th inner iteration, namely $$\rbf^l_k = ( (J^k)\tr J^k+\mu_k I)\ybf^{l} + \gbf^k,$$ for every $l\in\N$ we have the following
\begin{equation}\label{res1}\begin{aligned}
		&\|\rbf^{l+1}_k\| = \|( (J^k)\tr J^k+\mu_k I)\ybf^{l+1} + \gbf^k \| \\& = \|-(P^k+B^k+\mu_kI)( P^k+\mu_k I)^{-1}\left(\mathbf g^k+ B^k\ybf^l\right) + \gbf^k\| \\ &=  \|B^k\ybf^l + B^k(P^k+\mu_kI)^{-1}\lr{\gbf^k+B^k\ybf^l}\|\\
		&=\|B^k( P^k+\mu_kI)^{-1}\lr{(P^k+\mu_kI)\ybf^l+\gbf^k+B^k\ybf^l}\|\\
		&=\|B^k( P^k+\mu_kI)^{-1}\lr{((J^k)\tr J^k+\mu_kI)\ybf^l+\gbf^k}\|\\
		&\leq \|B^k( P^k+\mu_kI)^{-1}\|\|\rbf^l_k\| = \rho_k \|\rbf^l_k\|
	\end{aligned}
\end{equation} 
where we defined $\rho_k = \|B^k( P^k+\mu_kI)^{-1}\|.$

In the following, we use the fixed point iteration in \eqref{fixedpoint} to define an inexact LM method, PILM (\emph{Parallel Inexact Levenberg Marquardt Method}). More details regarding the implementation of the algorithm in the server/worker framework will be discussed in Section \ref{sec:PLS_results}, together with the numerical results. 
\newpage
\begin{alg}[PILM]\label{alg:PILM} $\ $\\
	\textbf{Parameters:} { $c>0, \{\ell_k\}_{k=0}^\infty\in\mathbb N,\ \{\varepsilon_k\}_{k=0}^\infty\in\R_{>0}$}\\
	\textbf{Iteration $k$:}
	\begin{algorithmic}[1]
		\STATE compute $ \Rbf^k$, $ J^k$, $P^k$, $B^k$, $ {  \mathbf g^k }$
		\STATE choose $\mu_k$
		\FOR{$i=1,\dots,K$}
		\STATE compute $\ybf^1_i$ such that $\left(P^k_i+\mu_kI_{n_i}\right)\ybf^1_i =-\mathbf g^k_i$
		\ENDFOR
		\FOR {$l=1,\dots, \ell_{k}-1$}
		\FOR{$i=1,\dots,K$}
		\STATE compute $\ybf^{l+1}_i$ such that $\left(P^k_i+\mu_kI_{n_i}\right)\ybf^{l+1}_i =-\lr{\mathbf g^k_i + \sum_{j=1}^KB_{ij}^{ k}\ybf^l_{ j}}$
		\ENDFOR
		\ENDFOR
		\STATE set $\dbf^k = (\ybf^{\ell_k}_1,\dots,\ybf^{\ell_k}_K)\tr$
		\STATE find the largest $\alpha_k\in\left\{2^{-j}\right\}_{j=0}^{\infty}$ such that
		\begin{equation}\label{alg:ArmijoIN}
			F  (\xbf^k+\alpha_k \dbf^k)\leq F  (\xbf^k)-c\alpha_k^2\|\mathbf g^k\|^2+\varepsilon_k 
		\end{equation}
		\STATE set $\xbf^{ k+1} = \xbf^k+\alpha_k \dbf^k$
	\end{algorithmic}
\end{alg}

\begin{remark}
	Let us consider the line search condition \eqref{alg:ArmijoIN}. For $\alpha_k$ that tends to zero, the term on the left-hand side tends to $F(\xbf^k)$, while the negative term in the right-hand side tends to zero. Since we assume that $\varepsilon_k>0$, one can always find $\alpha_k>0$ such that the line search condition \eqref{alg:ArmijoIN} is satisfied. Note that this argument holds even in case the direction $\dbf^k$ is not a descent direction for $F$ at $\xbf^k.$ In particular, Algorithm PILM is well defined.
\end{remark}
\begin{remark}
	The $K$ linear systems in line 4 are independent. In particular the rounds of the \emph{for} loop in lines 3-5 can be executed in a parallel fashion. The same holds also for the \emph{for} loop at lines 7-9. Additionally, the local chunks of $ J^k$, $P^k$ and $B^k$, { denoted by $ J^k_i, P_i^k $ and $ B_i^k$, } can be also computed in parallel, on each node $i=1,...,K$. This is demonstrated in Algorithm \ref{alg:PILM parallel}. The implementation is written in Python and uses the Message Passing Inteface (MPI) for parallelization. When building a parallel algorithm, it is important to consider the communication points, as they have a tendency to become performance bottlenecks. The proposed parallel algorithm has only a few necessary communication points: one at line 3, needed to compute $B$, one at line \REV{7}, needed to compute the sum and one at line 11, where $\dbf$ and $\gbf$ are gathered. The input data is read and distributed by the master process, while the parallel computations are performed equally on the master and working nodes. As the underlying communication approach is the master/worker framework, the communication is always between the master and the worker nodes, the worker nodes do not communicate mutually.
\end{remark}

\begin{alg}[PILM: Parallel implementation]\label{alg:PILM parallel} $\ $\\
	\textbf{Parameters:} { $c>0, \{\ell_k\}_{k=0}^\infty\in\mathbb N,\ \{\varepsilon_k\}_{k=0}^\infty\in\R_{>0}$}\\
	\textbf{Iteration $k$:}
	\begin{algorithmic}[1]
		\STATE compute $ \Rbf^k$
		\STATE on each node $i=1,...,K$:
		\begin{ALC@g}
			\STATE compute $ J^k_i$, $P^k_i$, $B^k_i$, $\REV{g^k_i}$
		\end{ALC@g}
		\STATE choose $\mu_k$
		\STATE on each node $i=1,...,K$:
		\begin{ALC@g}
			\STATE compute $\ybf^1_i$ such that $\left(P^k_i+\mu_kI_{n_i}\right)\ybf^1_i =-\mathbf g^k_i$
		\end{ALC@g}
		\FOR {$l=1,\dots, \ell_{k}-1$}
		\STATE on each node $i=1,...,K$:
		\begin{ALC@g}
			\STATE compute $\ybf^{l+1}_i$ such that $\left(P^k_i+\mu_kI_{n_i}\right)\ybf^{l+1}_i =-\lr{\mathbf g^k_i + \sum_{j=1}^KB_{ij}^k\ybf^l_{ j}}$
		\end{ALC@g}
		\ENDFOR
		\STATE gather $\dbf^k = (\ybf^{\ell_k}_1,\dots,\ybf^{\ell_k}_K)\tr$ and $\mathbf g^k = (\mathbf g_1^k,\dots,\mathbf g_K^k )$
		\STATE find the largest $\alpha_k\in\left\{2^{-j}\right\}_{j=0}^{\infty}$ such that
		\begin{equation}\label{alg:ArmijoIN}
			F  (\xbf^k+\alpha_k \dbf^k)\leq F  (\xbf^k)-c\alpha_k^2\|\mathbf g^k\|^2+\varepsilon_k 
		\end{equation}
		\STATE set $\xbf^{\REV{k+1}} = \xbf^k+\alpha_k \dbf^k$
	\end{algorithmic}
\end{alg}

\section{Convergence Analysis}\label{sec:PLS_convergence}
The following assumptions are regularity assumptions commonly used in LM methods.
\begin{assumption}{E2} \label{ass:LSD1} The vector of residuals $ \Rbf : \R^n \to \R^m $ is continuously differentiable.
\end{assumption}
\begin{assumption}{E3} \label{ass:LSD2} The Jacobian matrix $ J\in \R^{m\times n}$ of $ \Rbf $
	is $L$-Lipschitz continuous. That is, for every $\xbf,\mathbf y \in \mathbb{R}^n$
	$$\| J(\xbf) -  J(\mathbf{y}) \|\leq L\|\xbf - \mathbf{y}\|.$$
\end{assumption} 

\subsection{Global Convergence}

\begin{lemma}\label{lemma:PILM1}
	Assume that $\dbf^k$ is computed as in algorithm PILM for a given $\ell_k\in\N$. The following inequalities hold
	\begin{enumerate}[i)]
		\item for every $k\in\N_0$ $$\rho_k\leq\frac{\|B^k\|}{\mu_k}$$
		\item for every $k\in\N_0$
		$$\|\rbf_k^{\ell_k}\|\leq\rho_k^{\ell_k}\|\gbf^k\|$$
		\item for every $k\in\N_0$ 
		$$\|\dbf^k\|\leq\frac{(1+\rho_k^{\ell_k})}{\mu_k}\|\gbf^k\|$$
		\item for every $k\in\N_0$ 
		$$(\dbf^k)\tr\gbf^k\leq\lr{\frac{\rho_k^{\ell_k}}{\mu_k}-\frac{1}{\|J^k\|^2+\mu_k}}\|\gbf^k\|^2$$
		
		\item if $\mu_k$ in line 2 is chosen as 
		$\mu_k = \max\{\mu_{\min}, C_\mu\|{  B^k\|}\}$ for some $\mu_{\min}>0$ and $C_\mu>0,$ then \REV{
			$$\rho_k\leq \frac{1}{C_\mu}, $$
			and therefore
			$$\|\rbf_k^{\ell_k}\|\leq\lr{\frac{1}{C_\mu}}^{\ell_k}\|\gbf^k\|.$$}
		
	\end{enumerate}
	\begin{proof}
		By sub-multiplicativity of the norm, we have 
		\begin{equation}\label{res2}
			\rho_k = \|B^k(P^k+\mu_k I)^{-1}\|\leq \|B^k\|\|(P^k+\mu_k I)^{-1}\| \leq \frac{\|B^k\|}{\mu_k} 
		\end{equation}
		which is $i).$
		Using the definition of $\ybf^1$ in line { 6} of Algorithm { 2.2} we have
		\begin{equation}\label{res3}
			\begin{aligned}
				\|\rbf_k^1\| &= \|-(P^k+B^k+\mu_k I)(P^k+\mu_k I)^{-1} \gbf^k + \gbf^k\|= \\
				&= \|B^k(P^k+\mu_k I)^{-1}\|\|\gbf^k\| = \rho_k \|\gbf^k\|
			\end{aligned}
		\end{equation}
		This proves part $ii)$ of the thesis in case $\ell_k=1$. If $\ell_k>1$, recursively applying \eqref{res1}, and using the equality above, we get
		\begin{equation*}
			\begin{aligned}
				\|\rbf_k^{\ell_k}\| &\leq \rho_k\|\rbf_k^{\ell_k-1}\| \leq \rho_k^{\ell_k-1}\|\rbf_k^1\|  = \rho_k^{\ell_k}\|\gbf^k\|,
		\end{aligned}\end{equation*}
		and  $ii)$ is proved.
		
		By definition of $\dbf^k$ and $\rbf_k^{\ell_k}$ we have 
		\begin{equation}\label{res5}\dbf^k = ((J^k)\tr J^k+\mu_kI)^{-1}(-\gbf^k+\rbf_k^{\ell_k})\end{equation}
		Taking the norm, we have 
		\begin{equation*}
			\begin{aligned}\|\dbf^k\| &= \left\|((J^k)\tr J^k+\mu_k I)^{-1}\left(\gbf^k+\rbf_k^{\ell_k}\right)\right\| \\ 
				&\leq \|((J^k)\tr J^k+\mu_k I)^{-1}\|\left(\|\gbf^k\|+\|\rbf^{\ell_k}_k\|\right) \\
				&\leq \frac{(1+\rho_k^{\ell_k})}{\mu_k}\|\gbf^k\|,
			\end{aligned}
		\end{equation*}
		that is { $iii)$.}\\
		By \eqref{res5}, Cauchy-Schwartz inequality and $i)$, we have
		\begin{equation*}
			\begin{aligned}
				(\dbf^k)\tr\gbf^k &= (\gbf^k)\tr\left((J^k)\tr J^k+\mu_k\right)^{-1}\left(-\gbf^k+\rbf_k^{\ell_k}\right)  \\
				&\leq \lambda_{\max}\lr{(J^k)\tr J^k+\mu_k}^{-1}\|\REV{\rbf^{\ell_k}_k}\|\|\gbf^k\| \\ &-\lambda_{\min}\lr{(J^k)\tr J^k+\mu_k}^{-1}\|\gbf^k\|^2 \\
				&\leq\frac{1}{\mu_k}\rho_k^{\ell_k}\|\gbf^k\|^2-\frac{1}{\|J^k\|^2+\mu_k}\|\gbf^k\|^2 \\
				&\leq\lr{\frac{\rho_k^{\ell_k}}{\mu_k}-\frac{1}{\|J^k\|^2+\mu_k}}\|\gbf^k\|^2,
		\end{aligned}\end{equation*}
		and we have part { $iv)$} of the statement.\\
		To prove { $v)$} it is enough to notice that if $\mu_k\geq C_\mu\|B^k\|$ we have 
		$$\frac{\|B^k\|}{\mu_k} \leq \frac{1}{C_\mu},$$
		and thus { $v)$} follows directly from $i)$ and $ii)$.
	\end{proof}
\end{lemma}

\begin{theorem}
	Assume that Assumptions \ref{ass:LSD1} and \ref{ass:LSD2} hold, $\ell_k\geq \ell$ for every $k$, $\{\varepsilon_k\}$ is such that 
	$\sum_{k=0}^\infty \varepsilon_k< +\infty$, and that in line 2 $\mu_k$ is chosen as 
	$\mu_k = \max\{\mu_{\min}, C_\mu\|{ B^k}\|\}$ { with $ C_\mu >1. $ } Then for $\ell$ large enough we have that for every $\mathbf \xbf_0\in\R^n $, each accumulation point of the sequence $\{\xbf^k\}_{k=1}^\infty$ is a stationary point of ${ F}(\mathbf x).$
	\begin{proof}
		Applying recursively the line search condition \eqref{alg:ArmijoIN} we have that for every $k\in\N_0$
		\begin{equation}\label{INproof1}\begin{aligned}
				{ 0 \leq} F(\xbf^{k+1})&\leq F  (\xbf^k)-c\alpha_k^2\|\mathbf g^k\|^2+\varepsilon_k  \\ &\leq F(\xbf^0)-c\sum_{j=0}^k\alpha_j^2\|\gbf^j\|^2 + \sum_{j=0}^k\varepsilon_j.
			\end{aligned}
		\end{equation}
		Reordering inequality \eqref{INproof1} and taking the limit for $k\rightarrow+\infty$ we get
		$$\REV{c}\sum_{k=0}^{+\infty}\alpha_k^2\|\gbf^k\|^2\leq F(\xbf^0)+\sum_{k=0}^\infty \varepsilon_k< +\infty$$
		which implies that 
		$$\lim_{k\rightarrow+\infty}\alpha_k\|\gbf^k\| = 0.$$
		Let us consider $\xbf^*\in\Rn$ any accumulation point of the sequence $\{\xbf^k\}_{\REV{k=0}}^{\REV{\infty}}$. By definition of accumulation point, there exists an infinite subset of indices $\calK_0\subseteq\N_0$ such that the subsequence $\{\xbf^k\}_{k\in\calK_0}$ converges to $\xbf^*$. The limit above implies
		$$\lim_{k\in\calK_0}\alpha_k\|\gbf^k\| = 0.$$
		
		If there exists $\alpha>0$ such that $\alpha_k\geq\alpha$ for every index $k\in\calK_0$, then, by continuity of the gradient $\gbf$, 
		$$0=\lim_{k\in\calK_0}\alpha_k\|\gbf^k\| \geq \alpha\lim_{k\in\calK_0}\|\gbf^k\| = \alpha\|\gbf(\xbf^*)\| $$
		and therefore $\xbf^*$ is a stationary point of $F$.
		If such $\alpha>0$ does not exist, then one can find $\calK_1\subseteq\calK_0$ infinite set of indices such that $\lim_{k\in\calK_1}\alpha_k = 0$ and $\alpha_k<1$ for every $k\in\calK_1.$
		In particular this implies that for $\hat\alpha_k = 2\alpha_k$ condition \eqref{alg:ArmijoIN} does not hold. 
		That is, for every $k\in\calK_1$
		$$F(\xbf^k+\hat\alpha_k\dbf^k)>F(\xbf^k)-c\hat\alpha_k^2\|\gbf^k\|^2+\varepsilon_k.$$
		Since $\varepsilon_k\geq 0$ this implies 
		$$F(\xbf^k+\hat\alpha_k\dbf^k)-F(\xbf^k)> -c\hat\alpha_k^2\|\gbf^k\|^2,$$
		and applying the mean value theorem we have, for some $s_k\in[0,1]$
		\begin{equation}\label{INproof2}\begin{aligned}
				-c\hat\alpha_k\|\gbf^k\|^2 < \frac{1}{\hat\alpha_k}\lr{F(\xbf^k+\hat\alpha_k\dbf^k)-F(\xbf^k)} = \gbf(\xbf^k+s_k\hat\alpha^k\dbf^k)\tr\dbf^k.
			\end{aligned}
		\end{equation} 
		Let us now consider $\dbf^k.$ By part { $iii)$} of Lemma \ref{lemma:PILM1}, the definition of $\mu_k$, and the fact that $C_\mu>1$, we have 
		$$\|\dbf^k\|\leq \frac{1+\rho_k^{\ell_k}}{\mu_k}\|\gbf^k\|\leq\frac{2}{\mu_k}\max_{k\in\calK_1} \|\gbf^k\| $$
		where the maximum in the last term of the inequality exists because $\{\xbf^k\}_{k\in\calK_1}$is a compact subset of $\Rn$ and the gradient is continuous. Since the maximum is finite and $\mu_k\geq\mu_{\min}>0$ we have that $\{\dbf^k\}_{k\in\calK_1}$ is a bounded subsequence of $\Rn$ and therefore it has an accumulation point $\dbf^*$. That is
		$$\lim_{k\in\calK_2}\dbf^k = \dbf^*$$
		for some $\calK_2\subseteq\calK_1$ infinite subset.\\
		Since $\lim_{k\in\calK_2}\alpha_k = 0$, by definition of $\hat\alpha_k$ we have $\lim_{k\in\calK_2}\hat\alpha_k = 0$ and thus $\lim_{k\in\calK_2}\xbf^k+s_k\hat\alpha_k\dbf^k = \xbf^*$, which in turn implies
		$$\lim_{k\in\calK_2}\gbf(\xbf^k+s_k\hat\alpha_k\dbf^k)\tr\dbf^k=(\gbf^*)\tr\dbf^*.$$
		{ Taking the limit for $k\in\calK_2$ in \eqref{INproof2}}, we then get \begin{equation}\label{INproof3}
			0\leq \lim_{k\in\calK_2}(\dbf^k)\tr\gbf^k = (\dbf^*)\tr\gbf^*.
		\end{equation}
		On the other hand, by Lemma \ref{lemma:PILM1} \REV{$iv)$} and the definition of $\mu_k$ we have that 
		\begin{equation}\label{INproof4}
			\begin{aligned}
				(\dbf^*)\tr\gbf^* &= \lim_{k\in\calK_2}(\dbf^k)\tr\gbf^k\leq \lim_{k\in\calK_2}\lr{\frac{\rho_k^{\ell_k}}{\mu_k}-\frac{1}{\|J^k\|^2+\mu_k}}\|\gbf^k\|^2  \\ &\leq\lim_{k\in\calK_2} \lr{\frac{\rho_k^{\ell_k}}{\mu_{\min}}-\frac{1}{\displaystyle\max_{k\in\calK_2}\|J^k\|^2+\mu_{\max}}}\|\gbf^k\|^2.
			\end{aligned}
		\end{equation}
		where $\mu_{\max} =  { \{\mu_{\min}, C_{\mu} \max_{k\in\calK_2}\|B^k\|\}}$ and the maxima exist by compactness of $\{\xbf^k\}_{k\in\calK_2}$ and continuity of $J(\xbf)$ and $B(\xbf).$\\
		
		If $\lim_{k\in\calK_2}\|\gbf^k\| = 0$ then, by uniqueness of the limit, we have that $$\|\gbf^*\| = \lim_{k\in\calK_1}\|\gbf^k\| = \lim_{k\in\calK_2}\|\gbf^k\| = 0$$
		and therefore $\xbf^*$ is a stationary point of $F.$ Otherwise, we proceed by contradiction. If $\|\gbf^k\|$ does not vanish for $k\in\calK_2$, there exist $\gamma>0$ and $\calK_3\subseteq\calK_2$ infinite sequence such that $\|\gbf^k\|\geq\gamma$ for every $k\in\calK_3.$
		Fix { some $ \nu $ satisfying $ 0 < \nu < (\max_{k \in {\cal K}_2} \|J^k\|^2 + \mu_{\max})^{-1}$}.  From Lemma \ref{lemma:PILM1} and we have that $\rho_k^{\ell_k}\leq C_\mu^{-\ell_k}$. Since $C_\mu>1$ and $\ell_k\geq \ell$ for every $k$, one can find $\ell$ large enough such that
		$$\lr{\frac{\rho_k^{\ell_k}}{\mu_{\min}}-\frac{1}{\max_{k\in\calK_2}\|J^k\|^2+\mu_{\max}}}\leq - \nu.$$
		For this choice of $\ell$, from \eqref{INproof4} we have, for every $k\in\calK_3$
		$$(\dbf^*)\tr\gbf^*\leq -\nu\gamma^2.$$
		This contradicts \eqref{INproof3} and therefore concludes the proof.
	\end{proof}
\end{theorem}

\subsection{Local Convergence}

Let $S$ denote the set of all stationary points of $\| R (\xbf)\|^2$, namely $S = \{\xbf\in\R^{N}|  J(\xbf)\tr  R (\xbf) = 0\}.\ $ Consider a stationary point $ \xbf^{*}\in S $ and { an} open ball $ B_r:=B(\xbf^{*},r) $ with radius $ r \in(0,1) $ around it. 
From now on, given a point $\xbf\in\Rn$ we denote with $\bar\xbf$ a point in $S$ that minimizes the distance from $\xbf$. That is, 
$$\|\xbf-\bar\xbf\| = \min_{ z\in S}\|\xbf-\zbf\| = \dist(\xbf,S).$$

Since $B_r$ is a bounded subset of $\Rn$, and $\Rbf$ and $J$ are continuous functions on $\Rn$, we have that there exist $R_{\max},L_2\geq 0$ such that  for every $\xbf\in B_r$ $\|\Rbf(\xbf)\|\leq R_{\max}$ and $\|J(\xbf)\|\leq L_2.$
The following Lemma includes a set of inequalities, proved in \cite{santos:local} that are a direct consequence of assumptions \ref{ass:LSD1}, \ref{ass:LSD2} on the bounded subset $B_r.$
\begin{lemma}\label{lemma_PIN_loc0}\cite{santos:local} If Assumptions \ref{ass:LSD1}-\ref{ass:LSD2} hold, then for every $\xbf,\ybf\in B_r$
	\begin{enumerate}[i)]
		\item  $\|\Rbf(x)-\Rbf(y)-J(\ybf)(\xbf-\ybf)\|\leq \frac{1}{2}L\|\xbf-\ybf\|^2$
		\item $\|\Rbf(\xbf)-\Rbf(\ybf)\|\leq L_2\|\xbf-\ybf\|$
		\item $\|\gbf(\xbf)-\gbf(\ybf)\|\leq L_3\|\xbf-\ybf\|$ with $L_3 = L_2^2+L_1R_{\max}$
		\item denoting with $L_4 = \frac{1}{2}L_1L_2$ $$\|\gbf(\ybf)-\gbf(\xbf) - J(\xbf)\tr J(\xbf)(\ybf-\xbf)\|\leq L_4\|\xbf-\ybf\|^2 + \|(J(\xbf)- J(\ybf))\tr \Rbf(\ybf)\|$$  
		\item for every $\bar\zbf\in B_r\cap S$
		\begin{equation*}
			\begin{aligned}
				\|(J(\xbf)-J(\ybf))\tr\Rbf(\ybf)\|&\leq  L_1L_2\|\xbf-\bar\zbf\|\|\ybf-\bar\zbf\| + L_1L_2\|\ybf-\bar\zbf\|^2 \\
				&+\|J(\xbf)\tr\Rbf(\bar\zbf)\| + \|J(\ybf)\tr\Rbf(\bar\zbf)\|
			\end{aligned}
		\end{equation*}
	\end{enumerate}
\end{lemma}

In the rest of the section we make the following additional assumptions, which are standard for the local convergence analysis of Levenberg-Marquardt method. In particular, Assumption \ref{ass:LSD8}, referred to in the literature as local error bound condition, is typically used in place of the nonsingularity of the Jacobian matrix, while Assumption \ref{ass:LSD9} is common for non-zero residual problems.

\begin{assumption}{E4} \label{ass:LSD8}
	There exists $\omega>0$ such that for every $\xbf \in { B_r}$
	$$\omega\dist(\xbf, S)\leq \| J(\xbf)\tr  R (\xbf)\|$$
\end{assumption}
\begin{assumption}{E5} \label{ass:LSD9}
	There exist $\sigma>0$, \REV{$\delta\in[0,1]$} such that for every $\xbf \in { B_r} $ and every $\bar{\zbf} \in { B_r} \cap S$
	$$\|( J(\xbf) -  J(\bar{\zbf}))\tr  R (\bar{\zbf})\|\leq \sigma \|\xbf-\bar{\zbf}\|^{1+\delta}.$$
\end{assumption}
Notice that since $\bar\zbf$ is a stationary point of $\|R(\xbf)\|^2$, the inequality above is equivalent to $$\|J(\xbf)\tr  R (\bar{\zbf})\|\leq \sigma \|\xbf-\bar{\zbf}\|^{1+\delta}.$$
\begin{lemma} \label{lemma_PIN_loc1}Let us assume that \ref{ass:LSD1}-\ref{ass:LSD9} hold and that $\{\xbf^k\}_{\REV{k=0}}^{\REV{\infty}}$ is the sequence generated by Algorithm $PILM$ with $\alpha_k = 1$ for every $k$. Moreover, let us assume that $\xbf^k,\xbf^{k+1}\in B_r$ and $\|\dbf^k\|\leq c_1\dist(\xbf^k,S)$ for some constant $c_1\geq0$. Then the following inequality holds with $c_2 = L_4c_1^2+L_1L_2(1+c_1)+L_1L_2(1+c_1)^2$
	\begin{equation*}
		\begin{aligned}
			\omega\|\xbf^{k+1}-\bar\xbf^{k+1}\|&\leq c_2\|\xbf^k-\bar\xbf^k\|^2 + \left(\sigma+\sigma(1+c_1)^{1+\delta}\right)\|\xbf^k-\bar\xbf^k\|^{1+\delta}  \\ &+ (L_3\rho_k^{\ell_k} + c_1\mu_k)\|\xbf^k-\bar\xbf^k\|{ . }
		\end{aligned}
	\end{equation*}
	\begin{proof}
		By the triangular inequality, the assumptions of the Lemma, and the definition of $\bar\xbf^k$, we have the following inequalities
		\begin{equation*}
			\begin{aligned}
				&\|\xbf^{k+1}-\bar\xbf^k\|\leq \|\dbf^k\|+\|\xbf^k-\bar\xbf^k\| \leq (1+c_1)\|\xbf^k-\bar\xbf^k\|\\
				&\|(J^k)\tr\Rbf(\bar\xbf^k)\|\leq\sigma\|\xbf^k-\bar\xbf^k\|^{1+\delta}\\
				&\|(J^{k+1})\tr\Rbf(\bar\xbf^k)\|\leq\sigma\|\xbf^{k+1}-\bar\xbf^k\|^{1+\delta}\leq \sigma(1+c_1)^{1+\delta}\|\xbf^k-\bar\xbf^k\|^{1+\delta}.
			\end{aligned}
		\end{equation*}
		From part $iv)$ and $v)$ of Lemma \ref{lemma_PIN_loc0},  using the inequalities above, we have 
		\begin{equation}\label{PILM_proff1}
			\begin{aligned}
				&\|\gbf^{k+1}-\gbf^k - (J^k)\tr J^k(\xbf^{k+1}-\xbf^{k})\| \\ &\leq L_4\|\xbf^{k+1}-\xbf^k\|^2 + \|(J^k- J^{k+1})\tr \Rbf^{k+1}\|  \\ &\leq
				L_4\|\dbf^k\|^2 + L_1L_2\|\xbf^k-\bar\xbf^k\|\|\xbf^{k+1}-\bar\xbf^k\| + L_1L_2\|\xbf^{k+1}-\bar\xbf^k\|^2 \\
				&+\|(J^k)\tr\Rbf(\bar\xbf^k)\| + \|(J^{k+1})\tr\Rbf(\bar\xbf^k)\| \\
				&\leq c_2 \|\xbf^k-\bar\xbf^k\|^2 + \left(\sigma+\sigma(1+c_1)^{1+\delta}\right)\|\xbf^k-\bar\xbf^k\|^{1+\delta}.
			\end{aligned}
		\end{equation}
		
		From part $iii)$ in Lemma \ref{lemma_PIN_loc0}, using the fact that $$\left((J^k)\tr J^k+\mu_k I\right)\dbf^k = -\gbf^k+\REV{\rbf^{\ell_k}_k},$$ that $\bar\xbf^k$ is a stationary point of $F$, and the assumption over $\dbf^k$, we get
		\begin{equation}\label{PILM_proff2}
			\begin{aligned}
				&\|\gbf^k+(J^k)\tr J^k\dbf^k\|\leq \|\gbf^k+\left((J^k)\tr J^k+\mu_k I\right)\dbf^k\| + \mu_k \|\dbf^k\|  \\ & \leq  
				\|\REV{\rbf^{\ell_k}_k}\| + \mu_k \|\dbf^k\| \leq \rho_k^{\ell_k}\|\gbf^k\| + \mu_kc_1\|\xbf^k-\bar\xbf^k\|  \\
				&\leq \rho_k^{\ell_k}\|\gbf^k - \gbf(\bar\xbf^k)\| + \mu_kc_1\|\xbf^k-\bar\xbf^k\|\leq \left(L_3\rho_k^{\ell_k}+c_1\mu_k\right)\|\xbf^k-\bar\xbf^k\|.
			\end{aligned}
		\end{equation}
		By Assumption \ref{ass:LSD8}, adding and subtracting $\gbf^k +(J^k)\tr J^k\dbf^k$ we have 
		$$\omega\dist(\xbf^{k+1},S)\leq\|\gbf^{k+1}\| \leq \|\gbf^{k+1}-\gbf^k - (J^k)\tr J^k\dbf^k\|+\|\gbf^k+(J^k)\tr J^k\dbf^k\|.$$
		Replacing the two terms of the right-hand side with the bounds found in \eqref{PILM_proff1} and \eqref{PILM_proff2}, we get the thesis.
	\end{proof}
	
\end{lemma}

\begin{lemma}\label{lemma_PIN_loc2}
	If Assumptions \ref{ass:LSD1}-\ref{ass:LSD9} hold and $\{\xbf^k\}_{\REV{k=0}}^{\REV{\infty}}$ is the sequence generated by Algorithm PILM with $\alpha_k=1$, $\mu_k = \max\{\mu_{\min}, C_\mu\|B^k\|\}$, and $\ell_k\geq \ell$ for a given $\ell\in\N$ then, if $\{\xbf^k\}_{\REV{k=0}}^{\REV{\infty}}\subset B_r$, there exists $c_1>0$ such that for every iteration index $k$ 
	$$\|\dbf^k\|\leq c_1\dist(\xbf^k,S).$$
	\begin{proof}
		Let us denote with $p$ the rank of $J(\xbf^*)\tr J(\xbf^*)$ and let $\{\lambda^*_i\}_{i=1}^n = \eig(J(\xbf^*)\tr J(\xbf^*))$, in nonincreasing order. For a given iteration index $k$, let us consider the eigendecomposition of $(J^k)\tr J^k$ \begin{equation}\label{eigendec}
			(J^k)\tr J^k = \lr{Q^k_1,Q^k_2}\lr{\begin{matrix}\Lambda^k_1 & \\ & \Lambda^k_2\end{matrix}} \lr{Q^k_1,Q^k_2}\tr\end{equation} with 
		$\Lambda^k_1 = \dia(\lambda^k_1,\dots\lambda^k_p)\in\R^{p\times p}$ and $\Lambda^k_2 = \dia(\lambda^k_{p+1},\dots\lambda^k_n)\in\R^{(n-p)\times (n-p)}$, where $\{\lambda^k_i\}_{i=1}^n = \eig(J(\xbf^k)\tr J(\xbf^k))$ again in nonincreasing order, and $Q^k_1\in\R^{n\times p},\ Q^k_2\in\R^{n\times (n-p)}$ { are orthonormal matrices.} By continuity of $J(\xbf)$ and of the eigenvalues over the entries of the matrix, we have that for $r$ small enough $\min_{i=1:p}\lambda^k_i\geq\lambda^*_p/2.$\\
		
		By \eqref{eigendec} we have, for $i=1,2$
		$$(Q^k_i)\tr(\REV{\rbf^{\ell_k}_k}-\gbf^k)=(Q^k_i)\tr\lr{(J^k)\tr J^k+\mu_k I}\dbf^k = \lr{\Lambda^k_i+\mu_k I}(Q^k_i)\tr\dbf^k. $$
		For $i=1$, by definition of $\rho_k^{\ell_k}$, the fact that $\rho_k^{\ell_k}\leq\rho_k^{\ell}\leq C_\mu^{-\ell}$, and the bound on $\lambda^k_p$ we have 
		\begin{equation}\label{dirPILM0}
			\begin{aligned}
				&\|(Q^k_1)\tr \dbf^k\| { = } \|\lr{\Lambda^k_1+\mu_k I}^{-1}(Q^k_1)\tr(-\gbf^k+\REV{\rbf^{\ell_k}_k})\| \\ &\leq
				\frac{1}{\lambda_{\min}(\Lambda^k_1+\mu_k)}\|-\gbf^k+\REV{\rbf^{\ell_k}_k}\| \\ & \leq \frac{2}{\lambda^*_p}(1+\rho_k^{\ell_k})\|\gbf^k\|\leq \frac{4L_3\REV{C_{\mu}^{-\ell}}}{\lambda^*_p}\|\xbf^k - \bar\xbf^k\|.
			\end{aligned}
		\end{equation}
		For $i=2$, by { Lemma 3.1}, Lemma \ref{lemma_PIN_loc0}, Assumption \ref{ass:LSD9}, and the fact that $\|\lr{\Lambda^k_2+\mu_k I}^{-1}\Lambda^k_2\|\leq1$, we have
		\begin{equation}\label{dirPILM1}
			\begin{aligned}
				&\|(Q^k_2)\tr \dbf^k\|{ = }  \|\lr{\Lambda^k_2+\mu_k I}^{-1}(Q^k_2)\tr(-\gbf^k+\REV{\rbf^{\ell_k}_k})\| \\ 
				&\leq
				\|\lr{\Lambda^k_2+\mu_k I}^{-1}(Q^k_2)\tr(\gbf^k-\gbf(\bar\xbf^k)-(J^k)\tr J^k(\xbf^k-\bar\xbf^k))\|  \\ & +
				\|\lr{\Lambda^k_2+\mu_k I}^{-1}(Q^k_2)\tr(J^k)\tr J^k(\xbf^k-\bar\xbf^k)\| +  \frac{1}{\mu_k}\|\REV{\rbf^{\ell_k}_k}\|   \\  & \leq
				\frac{1}{\mu_k}\|\gbf^k-\gbf(\bar\xbf^k)-(J^k)\tr J^k(\xbf^k-\bar\xbf^k)\|  \\ & + \|\lr{\Lambda^k_2+\mu_k I}^{-1}\Lambda^k_2(Q^k_2)\tr(\xbf^k-\bar\xbf^k)\| + \frac{\rho_k^{\ell_k}}{\mu_k}\|\gbf^k\| \\ & \leq
				\frac{1}{\mu_k}\lr{L_4\|\xbf^k-\bar\xbf^k\|^2+\|(J^k)\tr\Rbf(\bar\xbf^k)\|}+\lr{1+\frac{\rho_k^{\ell_k}L_3}{\mu_k}}\|\xbf^k-\bar\xbf^k\| \\ & \leq
				\frac{L_4}{\mu_k}\|\xbf^k-\bar\xbf^k\|^2+\frac{\sigma}{\mu_k}\|\xbf^k-\bar\xbf^k\|^{1+\delta}+\lr{1+\frac{\rho_k^{\ell_k}L_3}{\mu_k}}\|\xbf^k-\bar\xbf^k\|.
			\end{aligned}
		\end{equation}
		By assumption we have $\mu_k\geq\mu_{\min}$ and $\rho_k^{\ell_k}\leq\rho_k^{\ell}\leq C_\mu^{-\ell}$. Therefore, proceeding in the previous chain of inequalities
		\begin{equation}\label{dirPILM2}
			\begin{aligned}
				&\|(Q^k_2)\tr \dbf^k\|\leq \lr{1+\frac{L_4+\sigma+C_\mu^{-\ell}L_3}{\mu_{\min}}}\|\xbf^k-\bar\xbf^k\|.
			\end{aligned}
		\end{equation}
		By the fact that $(Q^k_1,Q^k_2)$ is an orthonormal matrix, putting together  \eqref{dirPILM0} and \eqref{dirPILM2} we get
		\begin{equation*}
			\begin{aligned}
				&\|\dbf^k\|^2\leq\|(Q^k_1)\tr \dbf^k\|^2+\|(Q^k_2)\tr \dbf^k\|^2 \\ &\leq
				\lr{\frac{{ 4C_\mu^{-\ell}}L_3}{\lambda^*_p}}^2\|\xbf^k - \bar\xbf^k\|^2+\lr{1+\frac{L_4+\sigma+C_\mu^{-\ell}L_3}{\mu_{\min}}}^2\|\xbf^k-\bar\xbf^k\|^2.
			\end{aligned}
		\end{equation*}
		By definition of $\bar\xbf^k$, this implies the thesis with 
		$$c_1 = \lr{\lr{\frac{{ 4C_\mu^{-\ell}}L_3}{\lambda^*_p}}^2+\lr{1+\frac{L_4+\sigma+C_\mu^{-\ell}L_3}{\mu_{\min}}}^2}^{1/2}.$$
		
	\end{proof}
\end{lemma}

Before we state the following Lemma, we notice that by Assumption \ref{ass:LS} and the definition of $L_2$, we have that $\|B(\xbf)\|\leq C_B L_2^2$ for every $\xbf\in B_r$. In particular, if $\xbf^{k}\in B_r$ and $\mu_k = \max\{\mu_{\min}, C_\mu\|B^k\|\}$, 

then 
\begin{equation}\label{boundmu}
	\mu_k\leq \mu_{\max},\enspace \text{with}\enspace \mu_{\max} = \max\{\mu_{\min},C_\mu C_B L_2^2\}.
\end{equation}
\begin{lemma}\label{lemma_PIN_loc3}
	Let Assumptions \ref{ass:LS}-\ref{ass:LSD9} hold and let us denote with $\{\xbf^k\}_{\REV{k=0}}^{\REV{\infty}}$ the sequence generated by Algorithm PILM with $\alpha_k=1$, $\mu_k = \max\{\mu_{\min}, C_\mu\|B^k\|\}$ for $C_\mu>1$ and $\ell_k\geq \ell$ for a given $\ell\in\N$. Moreover, let { us} assume that there exists $\nu\in(0,1)$ such that $\omega\nu>c_4$ with 
	$c_4 = \sigma+\sigma(1+c_1)^2+L_3C_\mu^{-\ell} + c_1\REV{\mu_{\max}}$.
	If $\xbf^0\in B(\xbf^*,\varepsilon)$ with $$\varepsilon\leq\min\left\{\frac{\omega\nu-c_4}{c_2}, \frac{r(1-\nu)}{1+c_1-\nu}\right\}$$ then we have that for every $k\in\N_0$
	\begin{enumerate}[i)]
		\item $\xbf^{k+1}\in B_r$
		\item $\dist(\xbf^{k+1},S)\leq\nu\dist(\xbf^k,S)$
		\item $\dist(\xbf^{k+1},S)\leq\varepsilon.$
	\end{enumerate}

	\begin{proof}
		We proceed by induction over $k.$
		By Lemma \ref{lemma_PIN_loc2} and the bound on $\varepsilon$
		\begin{equation}\label{PILMA1}
			\begin{aligned}
				\|\xbf^{1}-\xbf^*\|&\leq \|\dbf^0\|+\|\xbf^0-\xbf^*\| \\&\leq (1+c_1)\dist(\xbf^0,S) \leq (1+c_1)\varepsilon\leq r,
			\end{aligned}
		\end{equation}
		{ w}hich is $i)$ for $k=0.$
		Since $\xbf^0\in B_r$ the bound \eqref{boundmu} holds.  
		By Lemma \ref{lemma_PIN_loc1}, and the fact that $\xbf^0\in B_\varepsilon$ we then have 
		\begin{equation}\label{PILMA2}
			\begin{aligned}
				&\omega\dist(\xbf^{1},S)\leq c_2\|\xbf^0-\bar\xbf^0\|^2 + \left(\sigma+\sigma(1+c_1)^{1+\delta}\right)\|\xbf^0-\bar\xbf^0\|^{1+\delta}  \\ &+ (L_3\rho_0^{\ell_0} + c_1\mu_0)\|\xbf^0-\bar\xbf^0\| \\
				&\leq \left(c_2\varepsilon+\sigma+\sigma(1+c_1)^{1+\delta}+L_3\rho_0^{\ell_0} + c_1\mu_0\right)\|\xbf^0-\bar\xbf^0\| \\
				&\leq \left(c_2\varepsilon+\sigma+\sigma(1+c_1)^2+L_3C_\mu^{-\ell} + c_1\REV{\mu_{\max}}\right)\dist(\xbf^0,S) \\
				&\leq \omega\nu\dist(\xbf^0,S)
			\end{aligned}
		\end{equation}
		and therefore $ii)$ holds for $k=0$.
		To prove $iii)$ is now enough to notice that since $\nu<1$ we have $$\dist(\xbf^{1},S)\leq\nu\dist(\xbf^{0},S)\leq\nu\varepsilon\leq\varepsilon.$$
		
		Let us now assume that for every $i=1,\dots,k$ we have $\xbf^{i}\in B_r$, $\dist(\xbf^{i},S)\leq\nu\dist(\xbf^{i-1},S)$ and $\dist(\xbf^{i},S)\leq\varepsilon.$ We want to prove that the same holds for $i=k+1.$
		From the definition of $\xbf^{k+1}$, the triangular inequality, the inductive assumptions and Lemma \ref{lemma_PIN_loc2} we have
		\begin{equation}\label{PILMA1_2}
			\begin{aligned}
				&\|\xbf^{k+1}-\xbf^*\|\leq \sum_{i=0}^k\|\dbf^i\|+\|\xbf^0-\xbf^*\| \\&\leq 
				c_1\sum_{i=0}^k\dist(\xbf^i,S)+\|\xbf^0-\xbf^*\| \\&\leq 
				c_1\sum_{i=0}^k\nu^i\dist(\xbf^0,S)+\|\xbf^0-\xbf^*\| \\&\leq 
				\left(1+c_1\sum_{i=0}^k\nu^i\right)\varepsilon \leq \lr{1+\frac{c_1}{1-\nu}}\varepsilon.
			\end{aligned}
		\end{equation}
		Since $\varepsilon\leq\frac{r(1-\nu)}{1+c_1-\nu}$ we have that the right-hand side is smaller than $r$ and therefore $\xbf^{k+1}\in B_r.$\\
		Proceeding as in \eqref{PILMA2}, 
		\begin{equation}\label{PILMA2_2}
			\begin{aligned}
				&\omega\dist(\xbf^{k+1},S)\leq c_2\|\xbf^k-\bar\xbf^k\|^2 + \left(\sigma+\sigma(1+c_1)^{1+\delta}\right)\|\xbf^k-\bar\xbf^k\|^{1+\delta}  \\ &+ (L_3\rho_k^{\ell_k} + c_1\mu_k)\|\xbf^k-\bar\xbf^k\| \\
				&\leq \left(c_2\varepsilon+\sigma+\sigma(1+c_1)^2+L_3C_\mu^{-\ell} + c_1\REV{\mu_{\max}}\right)\dist(\xbf^k,S) \\
				&\leq \omega\nu\dist(\xbf^k,S)
			\end{aligned}
		\end{equation}
		which implies $ii).$
		Since $\nu<1$, part $iii)$ of the thesis follows directly from $ii)$ and the fact that $\dist(\xbf^{k},S)\leq\varepsilon$.

	\end{proof}
\end{lemma}

\begin{theorem}\label{thm:PILM}
	If the same Assumptions of Lemma \ref{lemma_PIN_loc3} hold, then 
	$\dist(\xbf^k,S)\to 0 $ linearly and $\xbf^k\to \bar{\xbf}\in S\cap \REV{B_r}$.
	\begin{proof}
		By part $ii)$ of Lemma \ref{lemma_PIN_loc3} we have that for every iteration index $k$ $$\dist(\xbf^{k+1},S)\leq\nu\dist(\xbf^{k},S)$$
		since $\nu<1$ this implies that the sequence $\dist(\xbf^{k},S)$ converges linearly to 0.
		To prove the second part of the thesis, let us consider $l,s\in\N_0$ with $l\geq s$. From Lemma \ref{lemma_PIN_loc3} we have 
		$$\|\xbf^l-\xbf^s\|\leq\sum_{i=s}^{l-1}\|\dbf^i\|\leq c_1\varepsilon\sum_{i=s}^{l-1}\nu^i = c_1\varepsilon\frac{\nu^s-\nu^l}{1-\nu}$$
		which implies that $\{\xbf^k\}_{\REV{k=0}}^{\REV{\infty}}$ is a Cauchy sequence and therefore is convergent. By Lemma \ref{lemma_PIN_loc3} $i)$ and the fact that $\dist(\xbf^k,S)\to 0 $,  the limit point of the sequence has to be in $S\cap \REV{B_r}$, which concludes the proof.
	\end{proof}
\end{theorem}

We saw so far that the near-separability property of the problem influences the choice of the damping parameter $\mu_k$. In order to ensure convergence of the fixed-point method in lines { 6-10} of Algorithm PILM, one has to chose $\mu_k$ large enough, depending on the norm of the matrix $B(\xbf).$ However, for classical Levenberg-Marquardt method, in order to achieve local superlinear convergence, the sequence of damping parameters typically has to vanish \cite{santos:local}. In the \REV{remainder} of this section we show that under a stronger version of the near separability condition, the proposed method achieves superlinear and quadratic local convergence with assumptions, other than the near separability one, that are analogous to those of the classical LM.
\begin{assumption}{E6}\label{ass:LSD10}
	For every $x\in B_r$ we have that $\| B(\xbf)\|<\lambda_{\min}( P^k){ .}$
\end{assumption}

If this assumption holds, then $\rho_k = \| B^k( P^k+\mu_kI)^{-1}\|<1$ for every choice of $\mu_k$.
The following Theorem shows that, whenever the previous Assumption holds and $\REV{\delta\in(0,1]}$ in Assumption \ref{ass:LSD9}, one can find a suitable choice of the damping parameter $\mu_k$ that ensures local superlinear convergence, provided that the number of inner iterations $\ell_k$ is large enough.

\begin{lemma}\label{lemma_PIN_loc5}
	Let Assumptions \ref{ass:LSD1}-\ref{ass:LSD10} hold with $\REV{\delta\in(0,1]}$ { in} \ref{ass:LSD9}, and let us denote with $\{\xbf^k\}_{\REV{k=0}}^{\REV{\infty}}$ the sequence generated by Algorithm PILM with $\alpha_k=1$, and $\mu_k$ such that $$\mmu\|\xbf^k-\bar\xbf^k\|^\delta\leq\mu_k \leq \Mmu\|\xbf^k-\bar\xbf^k\|^\delta$$ for $0<\mmu\leq \Mmu$. If for every $k$ the number of inner iterations $\ell_k$ is such that $\rho_k^{\ell_k}\leq C_\eta\|\xbf^k-\bar\xbf^k\|^\delta$ for some $C_\eta\geq0$, then there exists $c_1\geq0$ such that  $$\|\dbf^k\|\leq c_1\dist(\xbf^k,S)$$
	for every $k\in\N_0$.
\end{lemma}

\begin{remark}
	From Assumption \ref{ass:LSD8} and part $iii)$ in Lemma \ref{lemma_PIN_loc0}, we have that 
	$\omega\|\xbf^k-\bar\xbf^k\|\leq\|\gbf^k\|\leq L_3\|\xbf^k-\bar\xbf^k\|$
	therefore, taking $\mu_k = \bar\mu\|\gbf^k\|^\delta$ for any $\bar\mu>0$ satisfies the assumption of the Lemma, with $\mmu = \omega\bar\mu$ and $\Mmu = L_3\bar\mu.$ Analogously, taking $\ell_k$ such that $\rho_k^{\ell_k}\leq \frac{C_\eta}{L_3^\delta}\|\gbf^k\|^\delta$ yields $\rho_k^{\ell_k}\leq C_\eta\|\xbf^k-\bar\xbf^k\|^\delta$.
\end{remark}
\begin{proof}
	The proof is analogous to that of Lemma \ref{lemma_PIN_loc2}. Inequalities \eqref{dirPILM0} and \eqref{dirPILM1} still hold, as they are independent of the choice of $\mu_k$ and $\ell_k.$ With the assumptions of the current Lemma, from \eqref{dirPILM1} we get
	\begin{equation}\label{dirPILM4}
		\begin{aligned}
			&\|(Q^k_2)\tr\dbf^k\|\leq\frac{L_4\|\xbf^k-\bar\xbf^k\|^2}{\mmu\|\xbf^k-\bar\xbf^k\|^\delta}+\frac{\sigma\|\xbf^k-\bar\xbf^k\|^{1+\delta}}{\mmu\|\xbf^k-\bar\xbf^k\|^\delta} \\ & +\lr{1+\frac{C_\eta L_3\|\xbf^k-\bar\xbf^k\|^\delta}{\mmu\|\xbf^k-\bar\xbf^k\|^\delta}}\|\xbf^k-\bar\xbf^k\| 
			\\& \leq    \lr{1+\frac{L_4+\sigma+C_\eta L_3}{\mmu}}\|\xbf^k-\bar\xbf^k\|.
		\end{aligned}
	\end{equation}
	Putting together \eqref{dirPILM0}, \eqref{dirPILM4} and the fact that $\|\dbf^k\|^2 = \|(Q^k_1)\tr\dbf^k\|^2+\|(Q^k_2)\tr\dbf^k\|^2$ we get the thesis with 
	$$c_1 = \lr{\lr{\frac{\REV{4}L_3}{\lambda^*_p}}^2+\lr{1+\frac{L_4+\sigma+C_\eta L_3}{\mmu}}^2}^{1/2}.$$
\end{proof}

\begin{lemma}\label{lemma_PIN_loc6}
	Let us assume that the same hypotheses of Lemma \ref{lemma_PIN_loc5} hold, and let us define $c_3 = c_2+\sigma+\sigma(1+c_1)^{1+\delta}+L_3C_\eta+c_1\Mmu$ and fix $\nu\in(0,1).$
	If $$\varepsilon\leq\min\left\{\lr{\frac{\omega\nu}{c_3}}^{1/\delta}, \frac{r(1-\nu)}{1+c_1-\nu}\right\}$$ and $\xbf^0\in B(\xbf^*,\varepsilon)$ we have that for every $k\in\N_0$
	\begin{enumerate}[i)]
		\item $\xbf^{k+1}\in B_r$
		\item $\dist(\xbf^{k+1},S)\leq\nu\dist(\xbf^k,S)$
		\item $\dist(\xbf^{k+1},S)\leq\varepsilon.$
	\end{enumerate}
	\begin{proof}
		The proof proceeds analogously to that of Lemma \ref{lemma_PIN_loc3}. Let us first consider the case $k=0.$ Part $i)$ of the thesis is proved as in \eqref{PILMA1}.
		Proceeding as in \eqref{PILMA2}, using the assumptions on $\mu_k$ and $\rho_k^{\ell_k}$ and the fact that $\xbf^0\in B_\varepsilon$ we have 
		\begin{equation}\label{PILMA6}
			\begin{aligned}
				&\omega\dist(\xbf^{1},S)\\ &\leq 
				c_2\|\xbf^0-\bar\xbf^0\|^2 +       \left(\sigma+\sigma(1+c_1)^{1+\delta} L_3 C_\eta + c_1 \Mmu\right)\|\xbf^0-\bar\xbf^0\|^{1+\delta} \\
				&\leq c_3 \|\xbf^0-\bar\xbf^0\|^{1+\delta}\leq c_3\varepsilon^\delta\dist(\xbf^0,S),
			\end{aligned}
		\end{equation}
		and by the bound on $\varepsilon$ we have that $ii)$ holds for $k=0$. Since $\nu<1$, $iii)$ follows immediately.
		Let us now assume that for every $i=1,\dots,k$ we have $\xbf^{i}\in B_r$, $\dist(\xbf^{i},S)\leq\nu\dist(\xbf^{i-1},S)$ and $\dist(\xbf^{i},S)\leq\varepsilon.$ We want to prove that the same holds for $i=k+1.$
		\REV{From \eqref{PILMA1_2}, using the fact that $\varepsilon\leq\frac{r(1-\nu)}{1+c_1-\nu}$ we have
			\begin{equation*}
				\begin{aligned}
					&\|\xbf^{k+1}-\xbf^*\|\leq \lr{1+\frac{c_1}{1-\nu}}\varepsilon\leq r
				\end{aligned}
			\end{equation*}
			and therefore $\xbf^{k+1}\in B_r.$}\\
		Proceeding as in \eqref{PILMA6}, \REV{using the assumptions on $\rho_k^{\ell_k}$ and $\mu_k$, we have}
		\begin{equation}\label{PILMA7}
			\begin{aligned}
				\omega\dist(\xbf^{k+1},S)&\leq c_2\|\xbf^k-\bar\xbf^k\|^2 + \left(\sigma+\sigma(1+c_1)^{1+\delta}\right)\|\xbf^{\REV{k}}-\bar\xbf^{\REV{k}}\|^{1+\delta}  \\ &+ \left(L_3\rho_k^{\ell_k} + c_1\mu_k\right)\|\xbf^k-\bar\xbf^k\| \\
				&\REV{\leq  c_2\|\xbf^k-\bar\xbf^k\|^2 + \left(\sigma+\sigma(1+c_1)^{1+\delta}\right)\|\xbf^{\REV{k}}-\bar\xbf^{\REV{k}}\|^{1+\delta}}  \\ &\REV{+ \left(L_3C_\eta\|\xbf^k-\bar\xbf^k\|^\delta + c_1\Mmu\|\xbf^k-\bar\xbf^k\|^\delta\right)\|\xbf^k-\bar\xbf^k\|}\\
				&\leq c_3 \|\xbf^k-\bar\xbf^k\|^{1+\delta}\leq c_3\varepsilon^\delta\dist(\xbf^k,S),
			\end{aligned}
		\end{equation}
		which implies $ii).$
		Since $\nu<1$, part $iii)$ of the thesis follows directly from $ii)$ and the fact that $\dist(\xbf^{k},S)\leq\varepsilon$.
	\end{proof}
\end{lemma}

\begin{theorem}
	If the same Assumptions of Lemma \ref{lemma_PIN_loc6} hold, then 
	$\dist(\xbf^k,S)\to 0$ and $\xbf^k\to \bar{\xbf}\in S\cap \REV{B_r}$. The convergence of $\dist(\xbf^k,S)$ is superlinear if $\delta\in(0,1)$ and quadratic if $\delta=1.$
	\begin{proof}
		Convergence of $\dist(\xbf^k,S)$ to zero follows directly from part $ii)$ of the previous Lemma. Moreover, proceeding as in \eqref{PILMA7}, we get
		$$\dist(\xbf^{k+1},S)\leq \frac{c_3}{\omega}\|\xbf^k-\bar\xbf^k\|^{1+\delta} = \frac{c_3}{\omega}\dist(\xbf^{k},S)^{1+\delta} $$
		which implies superlinear convergence for $\delta\in(0,1)$ and quadratic convergence for $\delta = 1.$ Convergence of the sequence $\{\xbf^k\}_{k=0}^\infty$ is proved as in Theorem \ref{thm:PILM}.
	\end{proof}
\end{theorem}

\section{Implementation and Numerical Results} \label{sec:PLS_results}
In this section we discuss the parallel implementation of the proposed method and we present the results of a set of numerical experiments carried out to investigate the performance of method and the influence of the parameter $K.$

\begin{figure}[h]
	\centering
	\begin{subfigure}[b]{0.29\textwidth}
		\centering
		\includegraphics[width = \textwidth]{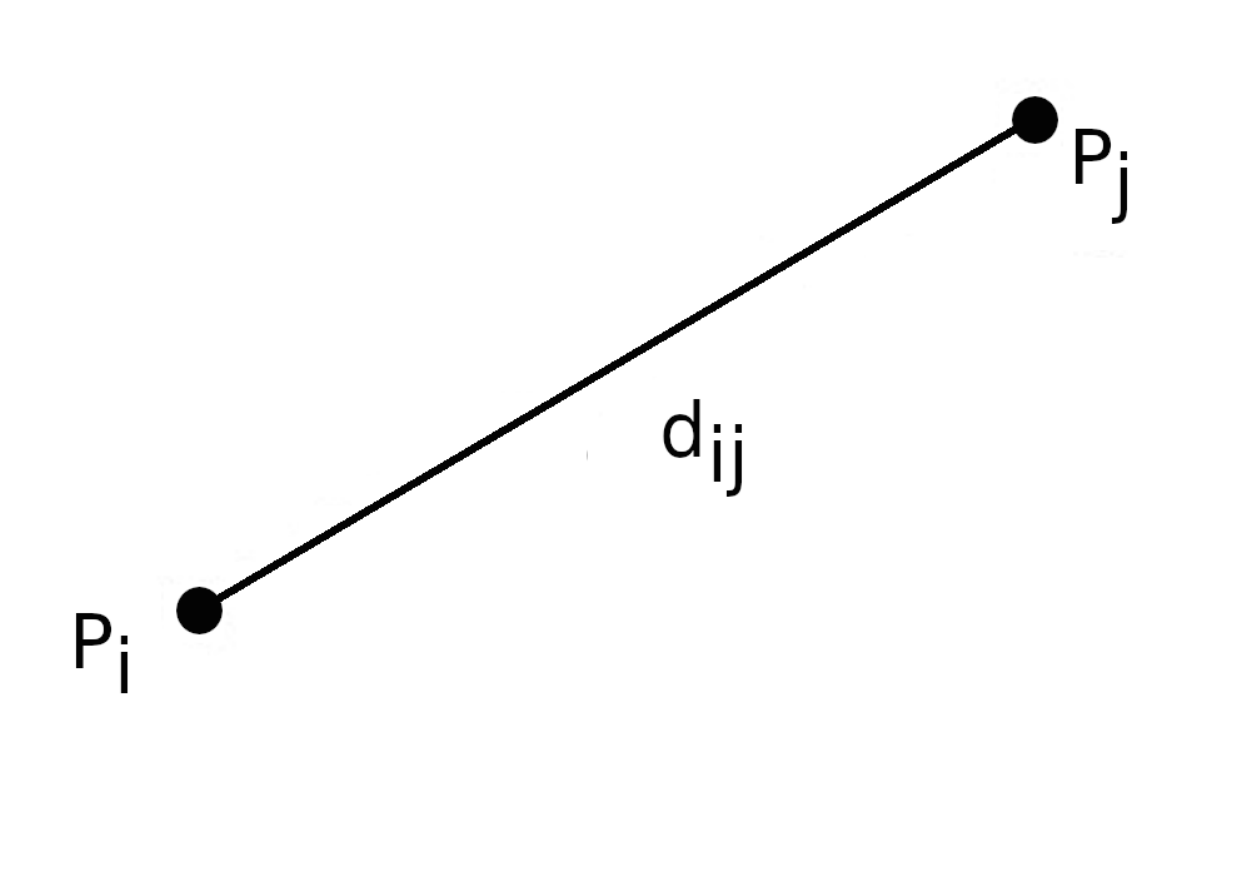}
		\caption{distance between $P_i$ and $P_j$\\}
	\end{subfigure}\hspace{0.03\textwidth}
	\begin{subfigure}{0.29\textwidth}
		\centering
		\includegraphics[width = \textwidth]{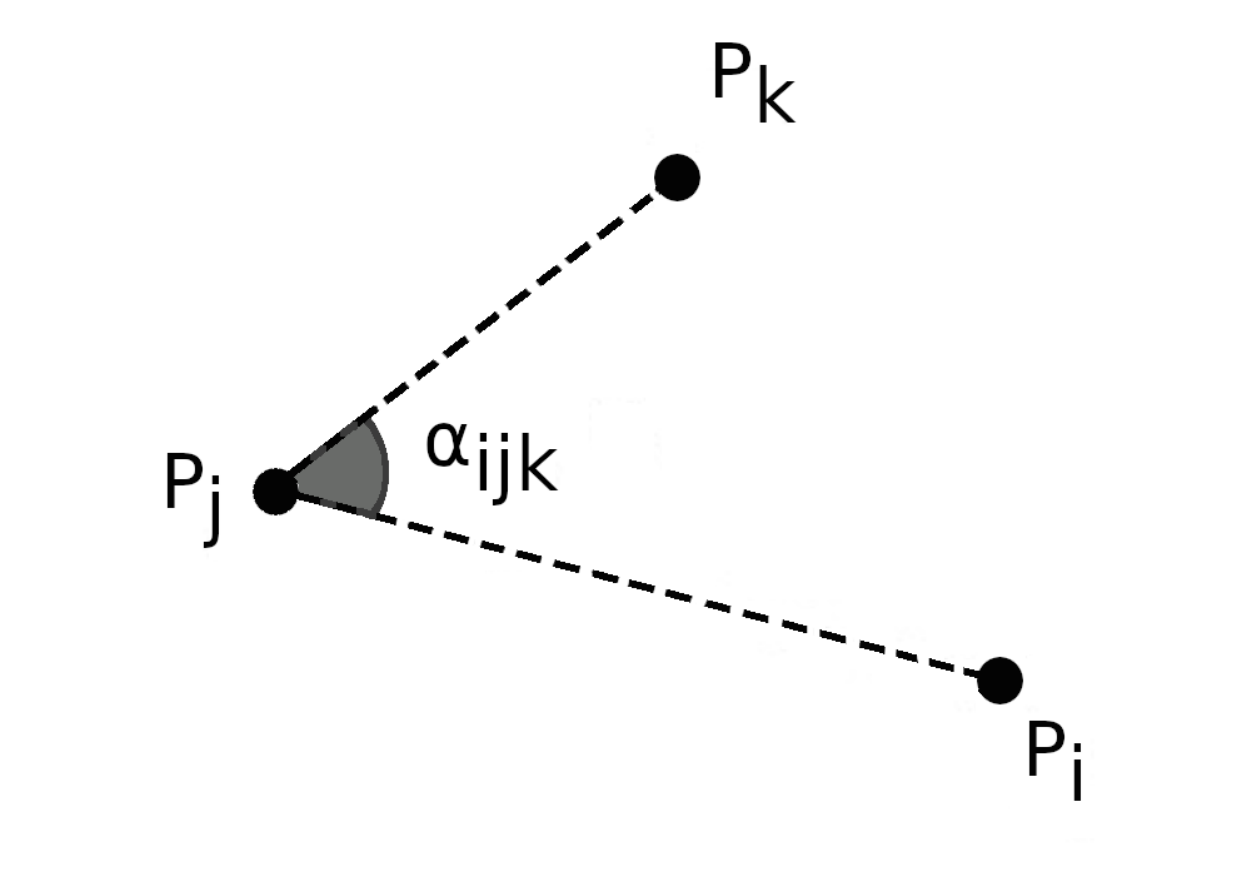}
		\caption{angle between by $P_i$ and $P_k$, centered at $P_j$\\}
	\end{subfigure}\hspace{0.05\textwidth}
	\begin{subfigure}{0.29\textwidth}
		\centering
		\includegraphics[width = \textwidth]{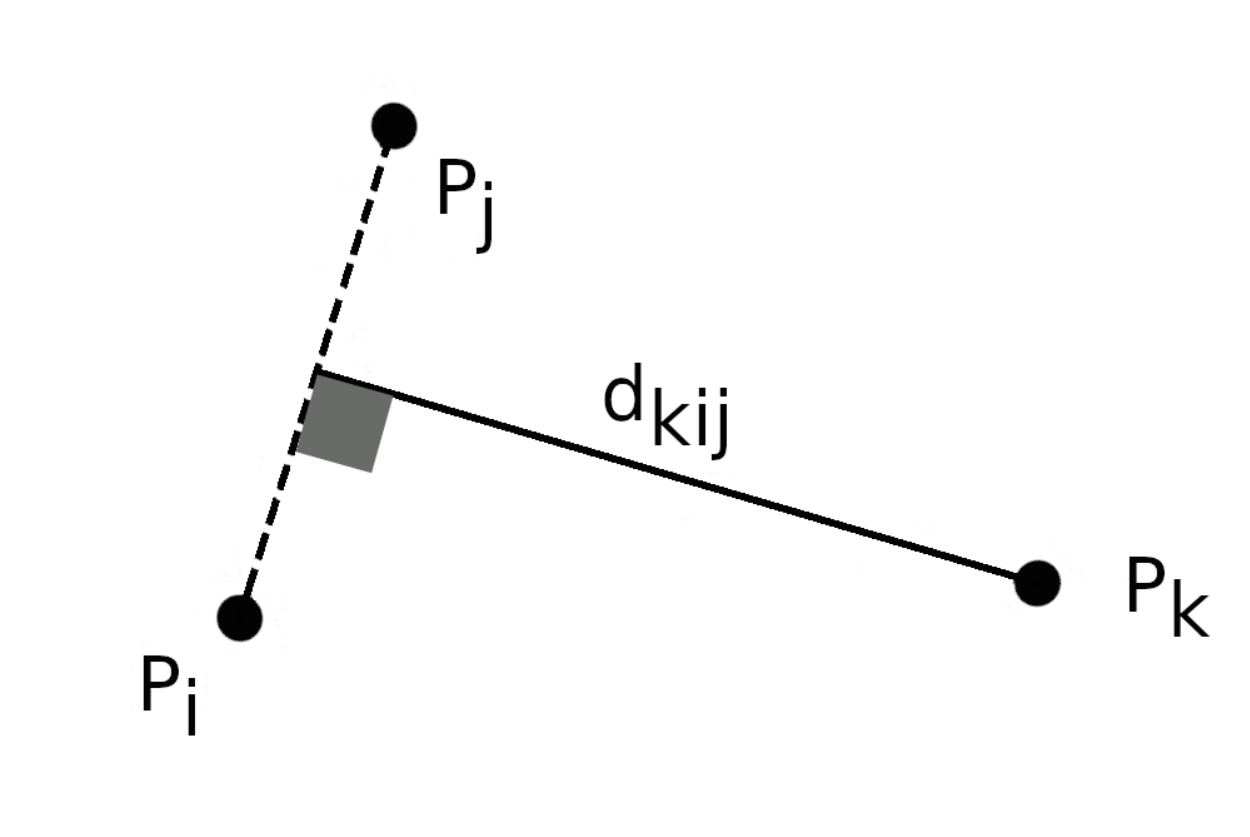}
		\caption{distance between $P_k$ and the line through $P_i$ and $P_j$ }
	\end{subfigure}
	
	\caption{Considered measurements.}\label{fig:measurements}
\end{figure}

We consider the least squares problems that arise from a Network Adjustment problem \cite{netadj}. Consider a set of points $\{P_1,\dots, P_n\}$ in $\mathbb{R}^2$ with unknown coordinates, and assume that a set of observations of geometrical quantities involving the points are available. Least Squares adjustments problems consist  {  into finding accurate coordinates of points using the available measurements,  by minimizing the residual with respect to the given observations in the least squares sense. } We consider here network adjustment problems with three kinds of observations: point-point distance, angle formed by three points and point-line distance, depicted in Figure \ref{fig:measurements}.

\REV{The problems we consider in the experiments are artificially generated and they simulate real problems of the same kind, in terms of average number of observations per point, average number of observation of each kind and near-separability. For additional details we refer to \cite{SLM}.}
Each problem is generated as follows, taking into account the information about average connectivity and structure of the network obtained from the analysis of real problems. Given the number of points $\hat n$ we consider a regular $2\sqrt{\hat n}\times 2\sqrt{\hat n}$ grid and we take $\{P_1,\dots, P_{\hat n}\}$ by uniformly sampling $25\%$ of the points on the grid. \REV{For each selected $P_i$ we define a neighborhood as follows: given an initial communication radius, we define the neighborhood as the set of points that are at a distance smaller than the communication radius. We repeat this procedure with increasing communication radius until the size of the neighborhood is at least 4. 
	Once the neighborhoods are defined, we generate a set of observations, proceeding as follows. We first select a point $P_i$, then we randomly select an observation kind. Specifically, we select point-point distance observations with probability 0.6, angle observations with probability 0.2 and point-line distance observations with probability 0.2. Once the kind is defined, we take the remaining points necessary to define the observation from the neighborhood of $P_i$. The observation is then generated by drawing a random value from a Gaussian distribution with mean equal to the true measurement (the actual euclidean distance between the two points, the actual angle formed by the three points, or the actual point-line distance) and standard deviation equal to 0.01 for distance observations and 1 for angle observations.}
We stop adding observations when every point is involved, on average, in 6 of them. We also generate \REV{coordinate observations for every point. That is, for each point $P_i$ we generate a random vector $(\hat x_i,\hat y_i)$ drawing $\hat x_i$ and $\hat y_i$  from the Gaussian distribution with mean equal to the true coordinate of the point and standard deviation equal to 1, for $99\%$ of the points, and 0.01 for the remaining $1\%.$}
We then consider the least square adjustment problem \cite{netadj} associated with the generated data.
That is, given the set of observations, the optimization problem is defined as a weighted least squares problem
\begin{equation}\label{LS_adj}
	\min_{\xbf\in\Rn}\frac{1}{2}\sum_{j=1}^m r_j(\xbf)^2 = \min_{\xbf\in\Rn} \frac{1}{2}\|\Rbf(\xbf)\|^2
\end{equation}
where $n=2\hat n$, $m$ is the number of observations, and $r_j(\xbf) = w_j^{-1}\widehat r_j(\xbf)$, with $\widehat r_j$ residual function of the $j$-th observation and $w_j$ corresponding standard deviation. 
For the three kinds of observations considered here (and represented in Figure \ref{fig:measurements}), the residuals $r_s$ are defined as follows:
\begin{itemize}
	\item distance between $P_i$ and $P_j$:\\
	\begin{equation}
		\label{eqn-measurement-1}
		r_s(x) = \sqrt{(x_{2i-1}-x_{2j-1})^2+(x_{2i}-x_{2j})^2} - d_{ij}
	\end{equation}
	\item angle between by $P_i$ and $P_k$, centered at $P_j$: 
	\begin{equation}
		\label{eqn-measurement-2}
		\begin{aligned}
			r_s(x) &= \atan((x_{2k}-x_{2j}), (x_{2k-1}-x_{2j-1}))+\\ &- \atan((x_{2i}-x_{2j}), (x_{2i-1}-x_{2j-1})) -\alpha_{ijk} 
		\end{aligned}
	\end{equation}
	\item distance between $P_k$ and the line through $P_i$ and $P_j$:
	\begin{equation}
		\label{eqn-measurement-3}
		r_s(x) = \frac{|(x_{2j-1}-x_{2i-1})(x_{2i}-x_{2k})-(x_{2i-1}-x_{2k-1})(x_{2j}-x_{2i})|}{\sqrt{(x_{2i-1}-x_{2j-1})^2+(x_{2i}-x_{2j})^2}} - d_{kij}.
	\end{equation}
\end{itemize}
Here, we assume that the variables in $\xbf$ are ordered in such a way that $(x_{2i-1}, x_{2i})$ corresponds to the coordinates of point $P_i$, for $i=1,\dots,\hat n,$
and quantities $d_{ij}$, $\alpha_{ijk}$ and 
$d_{kij}$ represent measurements illustrated in Figure \ref{fig:measurements}.\\

The proposed method is implemented in Python and all the tests are performed on the AXIOM computing facility consisting of 16 nodes (8 $\times$ Intel i7 5820k 3.3GHz and 8 $\times$ Intel i7 8700 3.2GHz CPU - 96 cores and 16GB DDR4 RAM/node) interconnected by a 10 Gbps network. \\

\begin{figure}[h]
	\centering
	\begin{subfigure}[b]{0.48\textwidth}
		\centering
		\includegraphics[width = 0.8\textwidth]{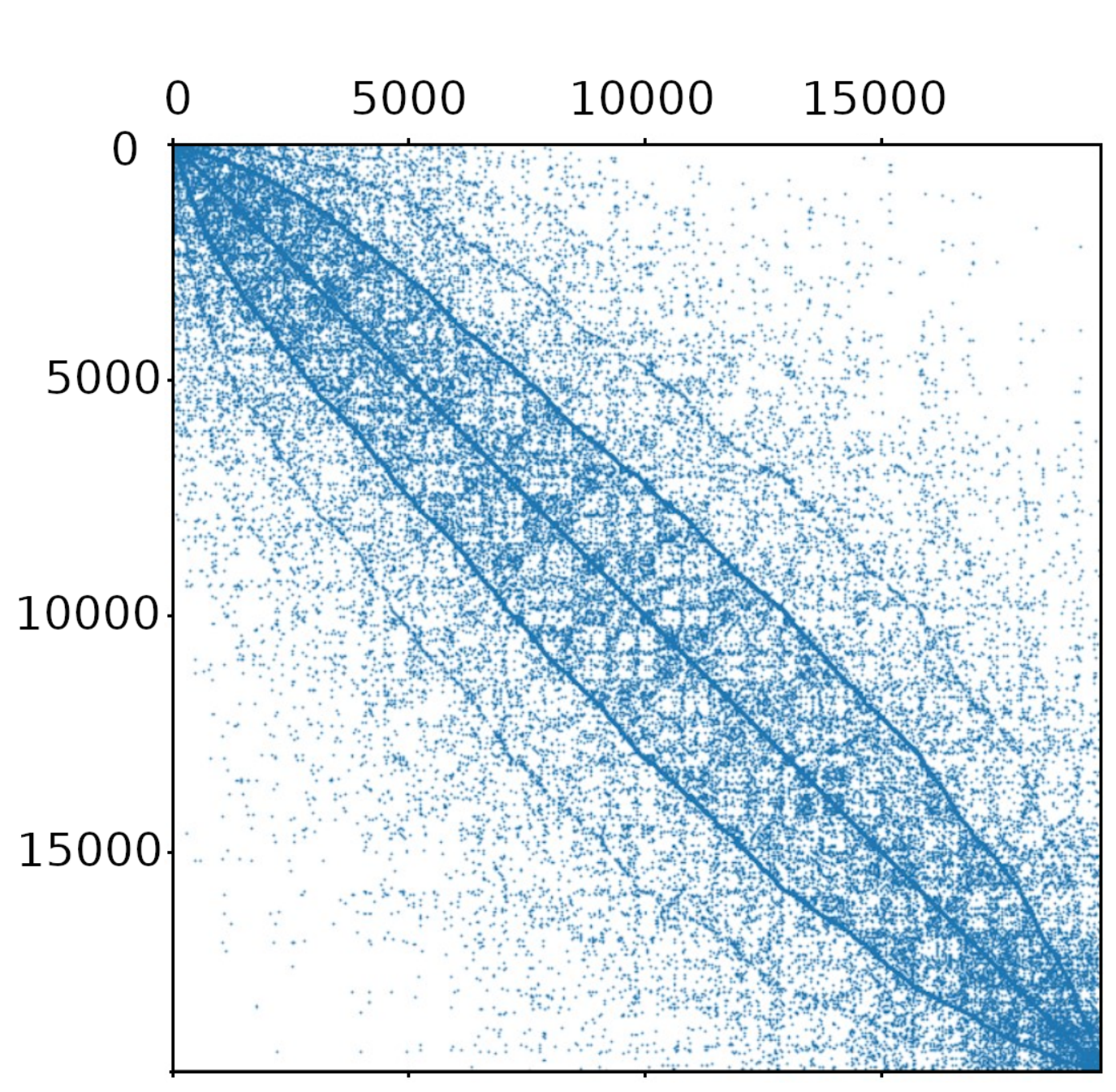}
		\caption{}
	\end{subfigure}
	\begin{subfigure}{0.48\textwidth}
		\centering
		\includegraphics[width = 0.8\textwidth]{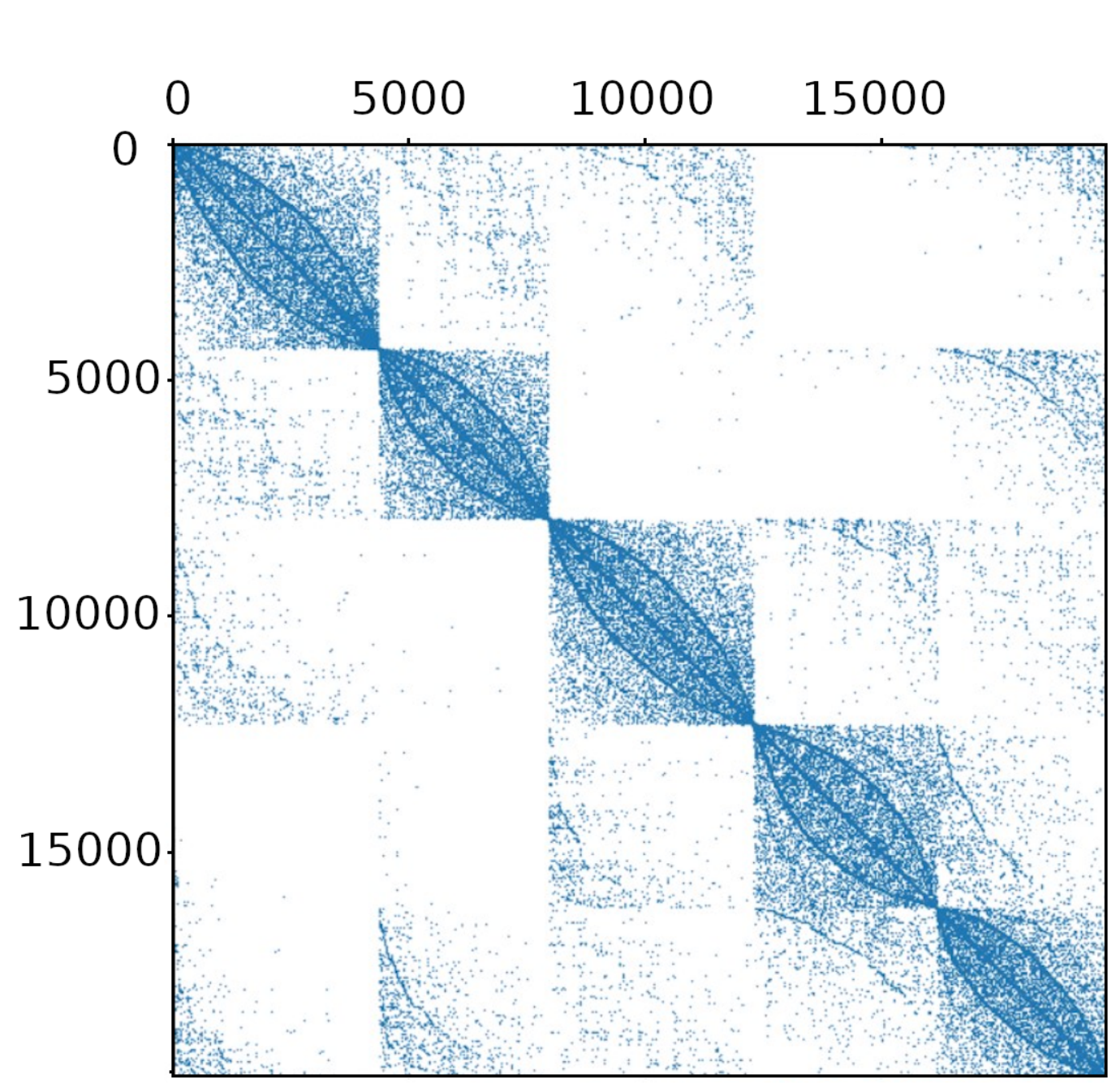}
		\caption{}
	\end{subfigure}
	
	\caption{Sparsity plot of the coefficient matrix of the system for one of the test problems.}\label{spymatrix}
\end{figure}

Algorithm \ref{alg:PILM} assumes that the number $K$ and the subsets $ I_s,E_s$, $s=1,\ldots,K$ of the variables and the residuals are given. In practice, the server computes the partition of the variables and defines the corresponding partition of the residuals as in \eqref{def_Es} and transmits them to the nodes. To compute the partition of the variables, in the tests that follow we use METIS \cite{metis} which is a graph-partitioning method that, given a network $\mathcal{G}$ and an integer $K$, separates the nodes into $K$ subsets in such a way that the cardinality of the subsets are similar, and that the number of edges between nodes in different subsets is approximately minimized. This method performs the splitting in a multi-level fashion starting from a coarse approximation of the graph and progressively refining the obtained partition, and it is therefore suitable for the partitioning of large networks, such as those that we consider. Moreover, the splitting of the variables and the residuals is computed only once at the beginning of the procedure. Overall, in the tests that we performed, the partitioning phase has limited influence over the total execution time. Nonetheless, the time employed by the algorithm to carry out the partitioning phase is included in the timings that we show below.\\
In Figure \ref{spymatrix} (a) we present the spyplot of the
matrix $ J\tr J$, for one of the considered test problems. In subfigure (b) we have the sparsity plot of the same matrix, where the rows and columns of the Jacobian have been reordered according to the results provided by METIS, with $K=5$. The block structure is that described in \eqref{JHB}.\\

In the following, given $i=1,\dots,K$ we denote with $\calN_i$ the set of indices $j\neq i$ such that there exists an observation in $\hat E$ involving variables in both $I_i$ and $I_j$. The Jacobian matrix and the derivatives of $F$ are computed as follows. For every $i=1,\dots,{ K}$ node $i$ computes $J_{i\Rbf_i}$ and $J_{i\rho}$ and shares $J_{i\rho}$ with the server, which then broadcasts $\{J_{i\rho}\}_{i=1}^K$ to the workers. Node $i$ then computes $g^k_i$, $P^k_i$ and $\{B^k_{ij}\}_{j\in\calN_i}$ according to \eqref{HandB_1}. Notice that $B^k_{ij}$ is nonzero only if $j\in\calN_j$. The local gradients $\gbf^k_i$ are transmitted by the workers to the server, that then broadcasts the aggregated gradient $\gbf^k.$
We observed that, compared to the approach where the server computes the whole Jacobian and then transmits it to the rest of the nodes, the distributed approach that we use leads to $50\%$ faster execution of this  phase.\\

To compute the direction $\dbf^k$ (lines 3-7 in Algorithm \ref{alg:PILM}) we proceed as follows. At the first inner iteration (line 3), node $i$ computes $\ybf_i^1$ solution of 
$$(P_i^k+\mu_kI_{n_i})\ybf_i^1 = -\gbf^k_i$$
which only involves quantities available to it. After solving this system, all nodes share their local solution with server. The server defines the aggregated vector $\ybf^l = (\ybf^l_1,\dots,\ybf^l_K)\tr$ and broadcasts it to the nodes. For all other inner iterations (line 8) each node $i$ first computes the right hand side 
$$(\gbf^k+B^k\ybf^l)_i = \gbf^k_i + \sum_{j\in\calN_i}B^k_{ij}\ybf^l_j$$
using the aggregated vector received from the server, then computes the new local estimate $\ybf_i^{l+1}$ as the solution of 
$$(P_i^k+\mu_kI_{n_i})\ybf_i^{l+1} = -(\gbf^k+B^k\ybf^l)_i.$$
Each node then sends the local vector to the server, which defines and shares the aggregated vector, and a new inner iteration begins. All the linear systems are solved with PyPardiso \cite{pypardiso}.\\

For all the communication phases we considered three approaches. The one mentioned above where the server broadcasts the aggregated quantities to all the nodes, the case where it sends to each node only the blocks that are necessary to them to perform their local computations, and the case where we define communicator between the workers, in such a way that node $i$ can share relevant quantities directly to the nodes in $\calN_i$. While in the second option the amount of exchanged data is smaller, we observed that the broadcasting approach results in practice in a significantly shorter communication time. Overall, the performance of the node-to-node approach was very similar to that of broadcasting and therefore we chose to continue with the broadcasting strategy, which is simpler from the point of view of the implementation.

\begin{figure}[h]
	\centering
	\includegraphics[width=0.5\linewidth]{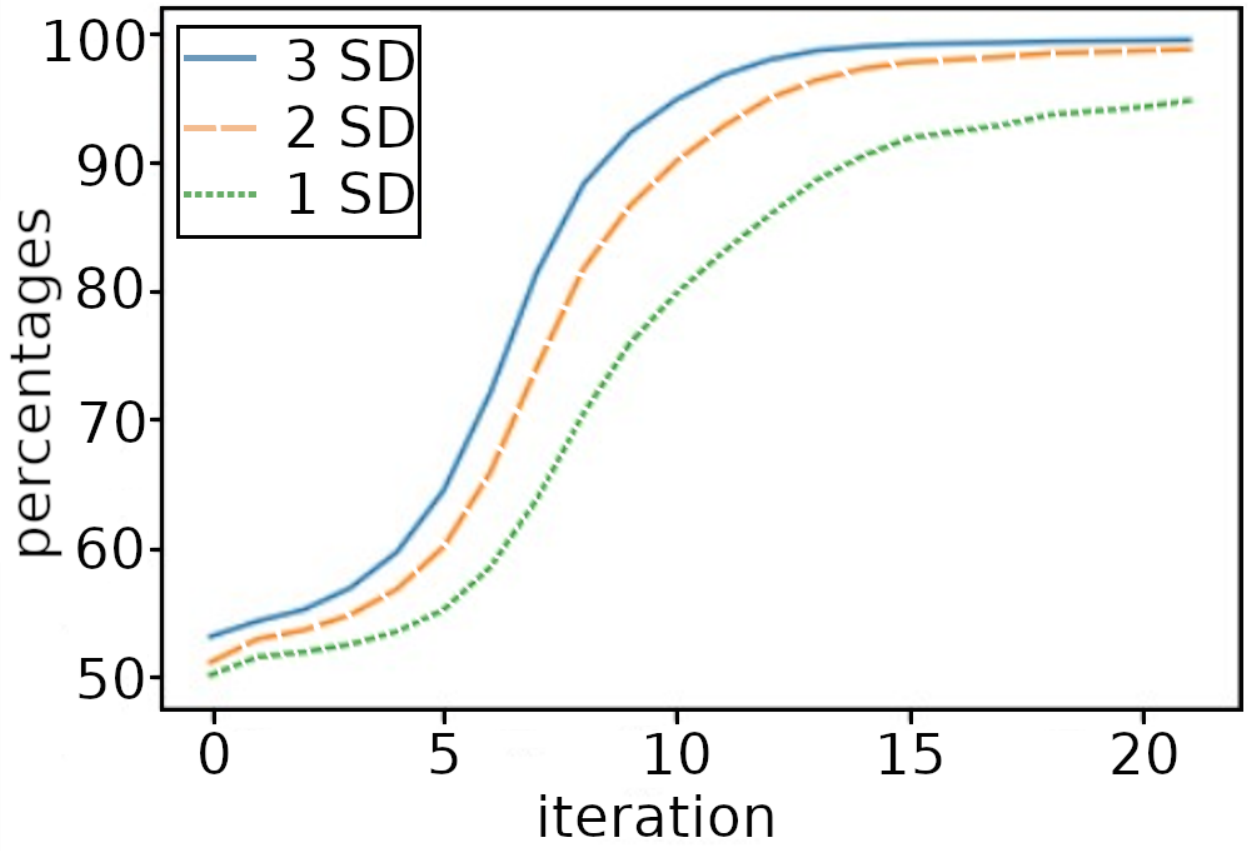}
	\caption{Percentage of residuals within 1, 2 and 3 standard deviations. Values of the percentages at each iteration }
	\label{fig:PILM_perc}
\end{figure}

\begin{figure}[h]
	\centering
	\includegraphics[width=0.8\linewidth]{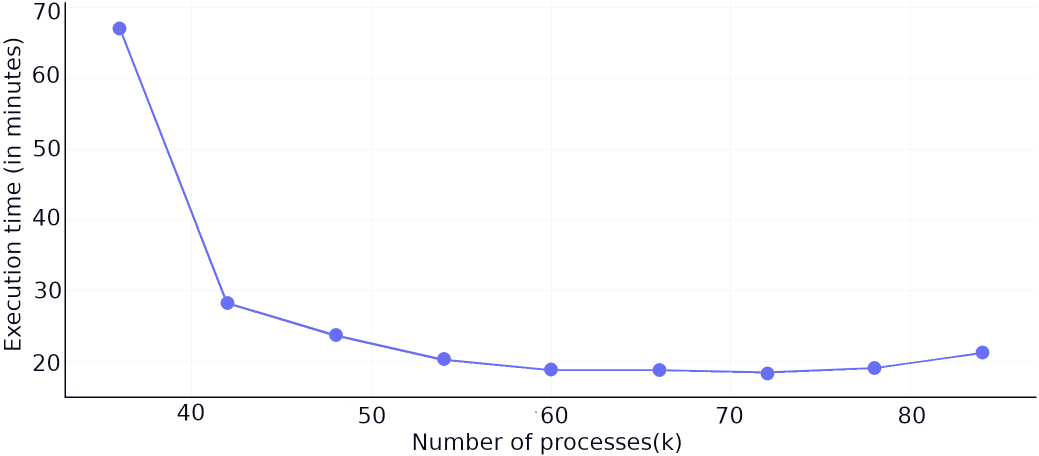}
	\caption{Execution time for different number of processors} 
	\label{fig:PILM1}
\end{figure}

We consider a network adjustment problem with $n=10^6$ and $m = 2.5\times10^6$ generated as described above, and we solve the problem for different values of the parameter $K$. { The initial guess is defined as the coordinate observations available in the problem description. The stopping criterion is as follows. For every $j=1,\dots,m$ we compare the $j$-th residual $r_j(\xbf^k)$ with the standard deviation $w_j$ associated with the $j$-th observation. The execution of the algorithm is then terminated when at least $68,\ 95\%$ and $99.5\% $ of the residuals is smaller than $w_j,\ 2w_j$ and $3w_j$, respectively. } To understand the behavior of the method, in Figure \ref{fig:PILM_perc} we plot the values of the three percentages above at each iteration, for K=60. In Figure \ref{fig:PILM1} we plot the execution time to arrive a termination for $K\in[35,85].$ The damping parameter $\mu_k$ is initialized as $10^5,$ which is the same order of magnitude as $\|\Rbf\|$. At each iteration we take $\mu_{k+1} = \mu_k/2$ if the accepted step size is larger than 0.5 and  $\mu_{k+1} = 2\mu_k$ otherwise, with safeguards $\mu_{\min} = 10^{-10}$ and $\mu_{\max} = 10^{10}$. The number of inner iteration is fixed to $\ell_k = 5$ for every $k$.

We can see that, starting from the smaller values of $K$, the execution time of the method decreases while $K$ increase, until reaching a plateau, after which it begins to increase. The plot shows the good performance of the proposed method. Smaller values of $K$ are omitted from the plot as the time necessary to arrive at termination becomes too large. In particular for $K=1$, which correspond to the centralized method, the execution time is orders of magnitude larger than for the values of $K$ included in the plot, and hence not comparable.\\
There are two main reasons behind the increase for large values of $K$. The first is that for larger values of $K$ the norm of the matrix $B$ is larger and therefore the fixed point method converges more slowly to the solution of the LM system. Since we are running a fixed number of inner iterations $\ell$ that does not depend on the number of nodes $K$, large values of $K$ result in a direction $\dbf^k$ that is a worse approximation of the LM direction and therefore the number of outer iterations needed by the method is larger. That is, after a certain point the overall computational cost increases because the saving induced by the fact that the linear systems solved by each node are smaller is not enough to balance the additional number of outer iterations. The second reason is common to all parallel methods: increasing the number of nodes $K$ increases the communication traffic and, when $K$ is too large, the time necessary to handle the additional communication overcomes the saving in terms of computation.\\ The fact that there is a plateau is also relevant from the practical point of view. The optimal value $K$ depends on the size $n$ but also on sparsity and the separability of the problem, and thus it may be hard to predict. However, the results show that the obtained timings on the considered problems are similar and nearly-optimal for a wide range of values of $K$, suggesting that an accurate choice of the number of nodes could in general not be necessary in order for the method to achieve a good performance.\\ Notice that while the choice of $\mu_k$ does not ensure theoretically $\rho_k<1$ at all iterations, it gives good results in practice, and does not require the computation of $\|B_k\|$, which may be expensive in the distributed framework.

As a comparison, Algorithm \ref{alg:PILM} was also implemented and tested on the same problem in a sequential fashion. That is, with only one machine performing the tasks for $i=1,\dots,K$ in sequence. The resulting timings were 228, 61.5 and 59.6 minutes for $K=45,80, 100$ respectively. Since these timings decrease for increasing $K$, this shows that the proposed method is effective, compared to classical LM method (equivalent to $K=1$), even when a parallel implementation is not possible in practice. Moreover, they show that the saving in time induced by the parallelization of the computation is significantly larger than the time necessary to handle the communication. \REV{For completeness, we report that the timings in the parallel case where 23.5 and 19.10 minutes for $K=42 $ and 78, respectively.}

\begin{figure}[h]
	\centering
	\begin{subfigure}[b]{0.48\textwidth}
		\centering
		\includegraphics[width = 1.1\textwidth]{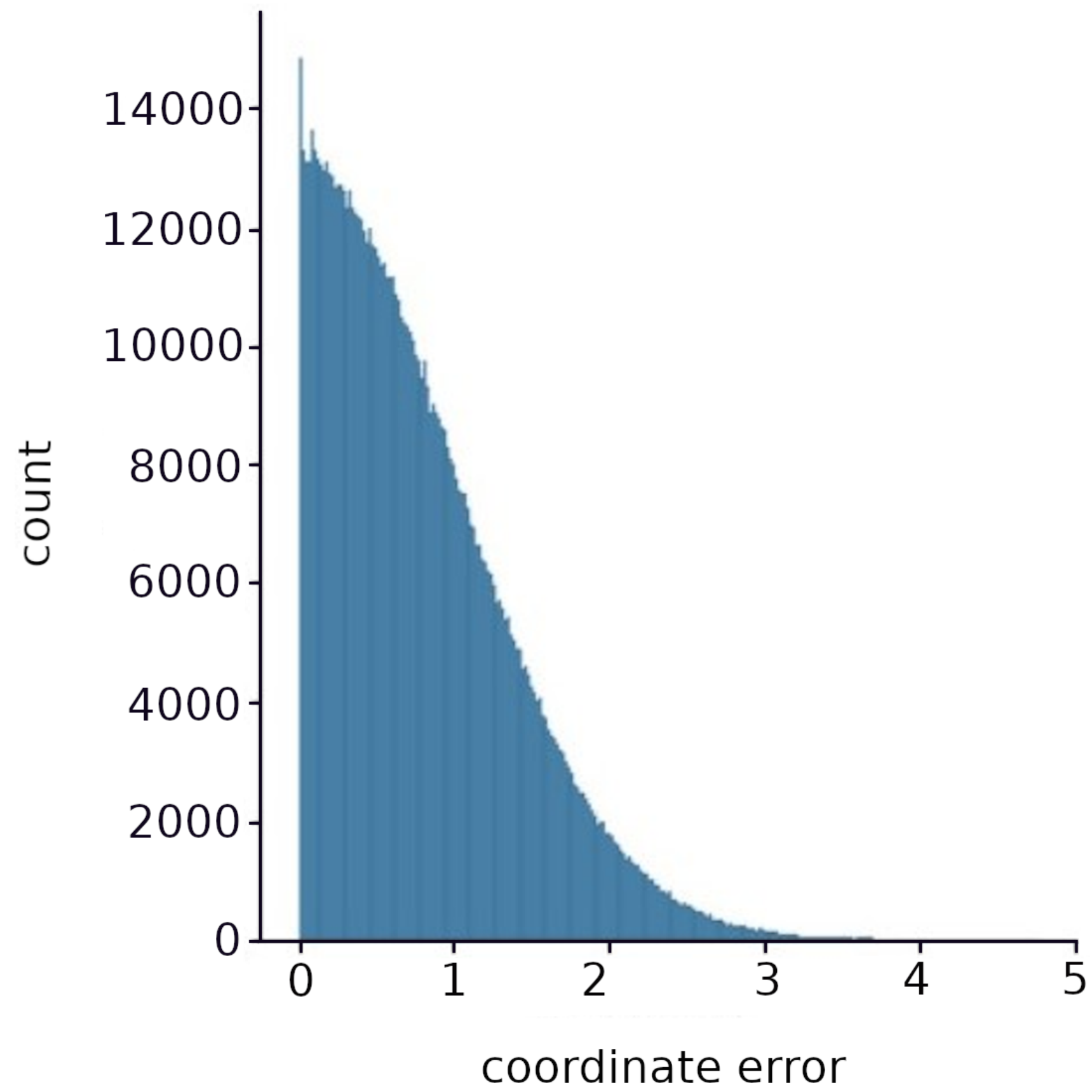}
	\end{subfigure}
	\begin{subfigure}{0.48\textwidth}
		\centering
		\includegraphics[width = 1.1\textwidth]{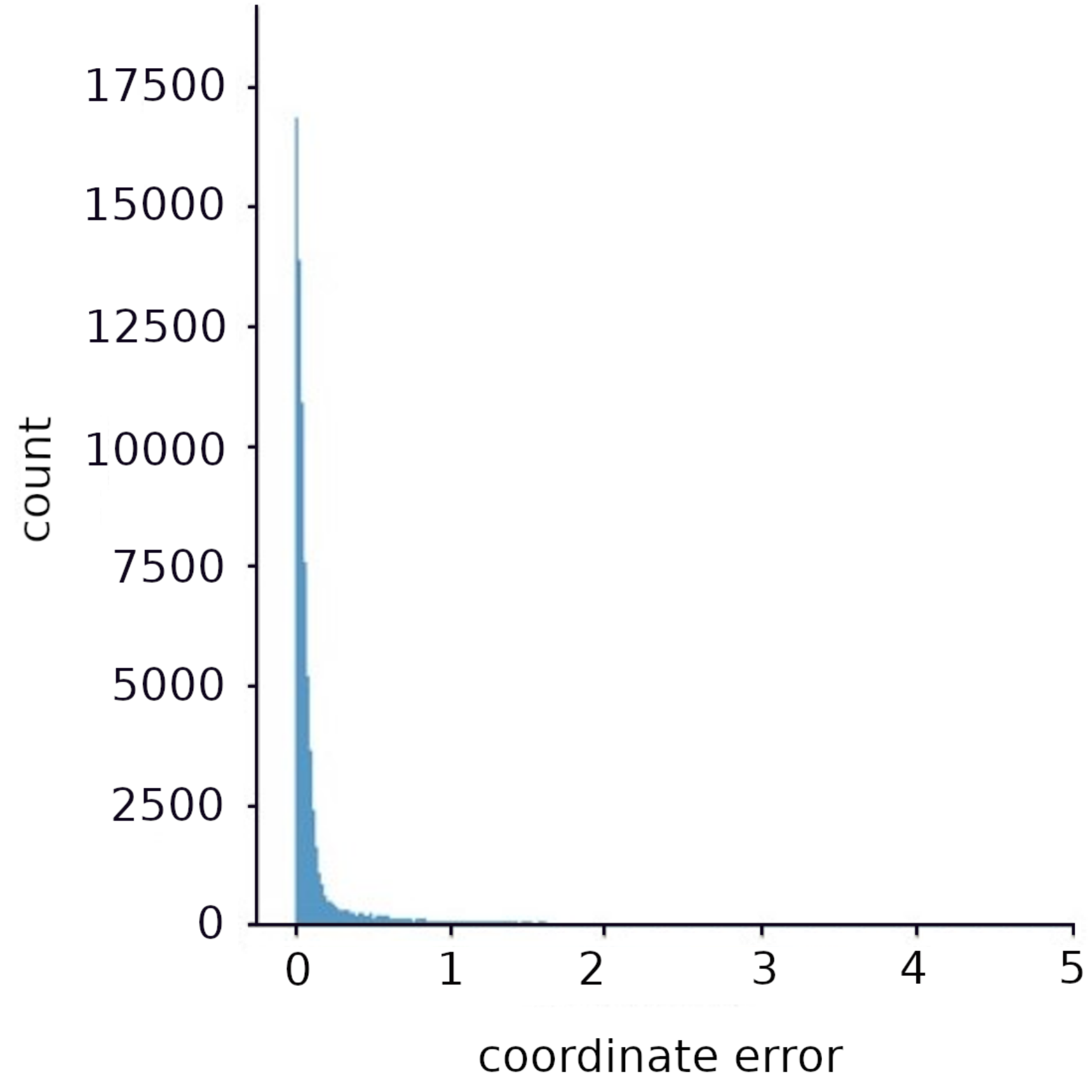}
	\end{subfigure}
	
	\caption{Distribution plot of the coordinate error at the initial guess (left) and at termination (right).}\label{err_500}
\end{figure}

\begin{figure}[h]
	\centering
	\begin{subfigure}[b]{0.327\textwidth}
		\centering
		\includegraphics[width = \textwidth]{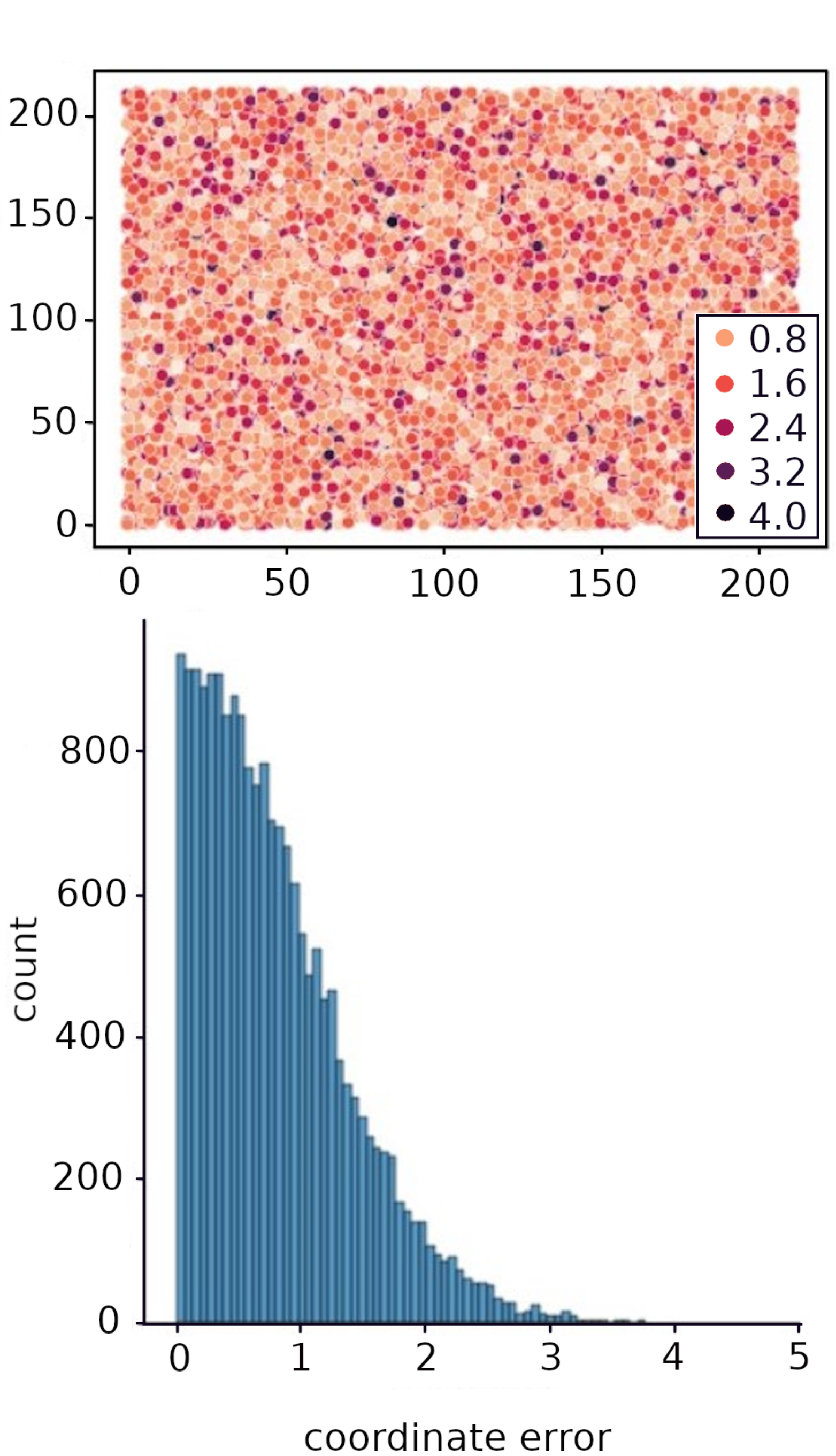}
		\caption{initial guess\\}
	\end{subfigure}
	\begin{subfigure}{0.327\textwidth}
		\centering
		\includegraphics[width = \textwidth]{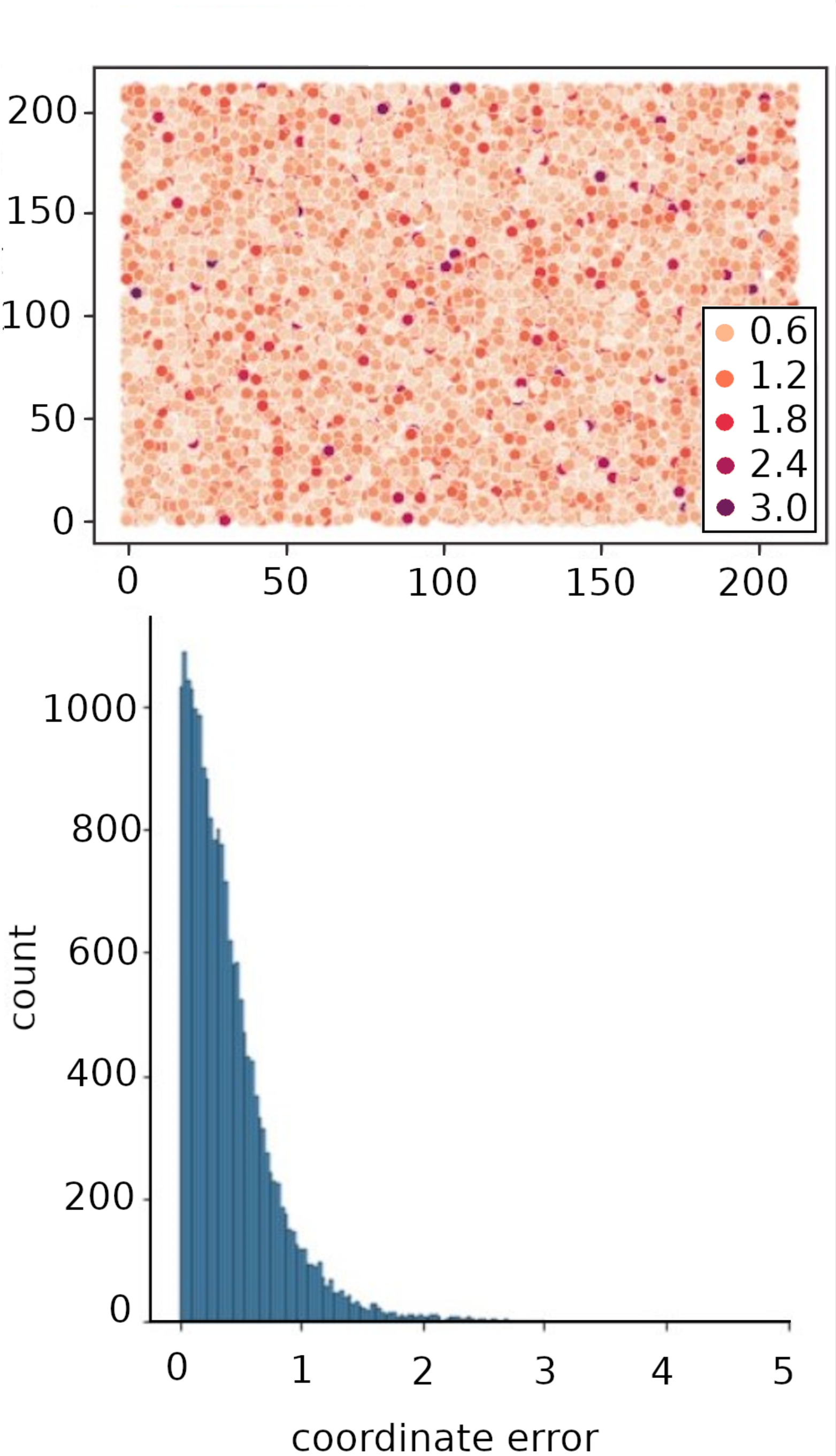}
		\caption{classical LM - 11.7s}
	\end{subfigure}
	\begin{subfigure}{0.327\textwidth}
		\centering
		\includegraphics[width = \textwidth]{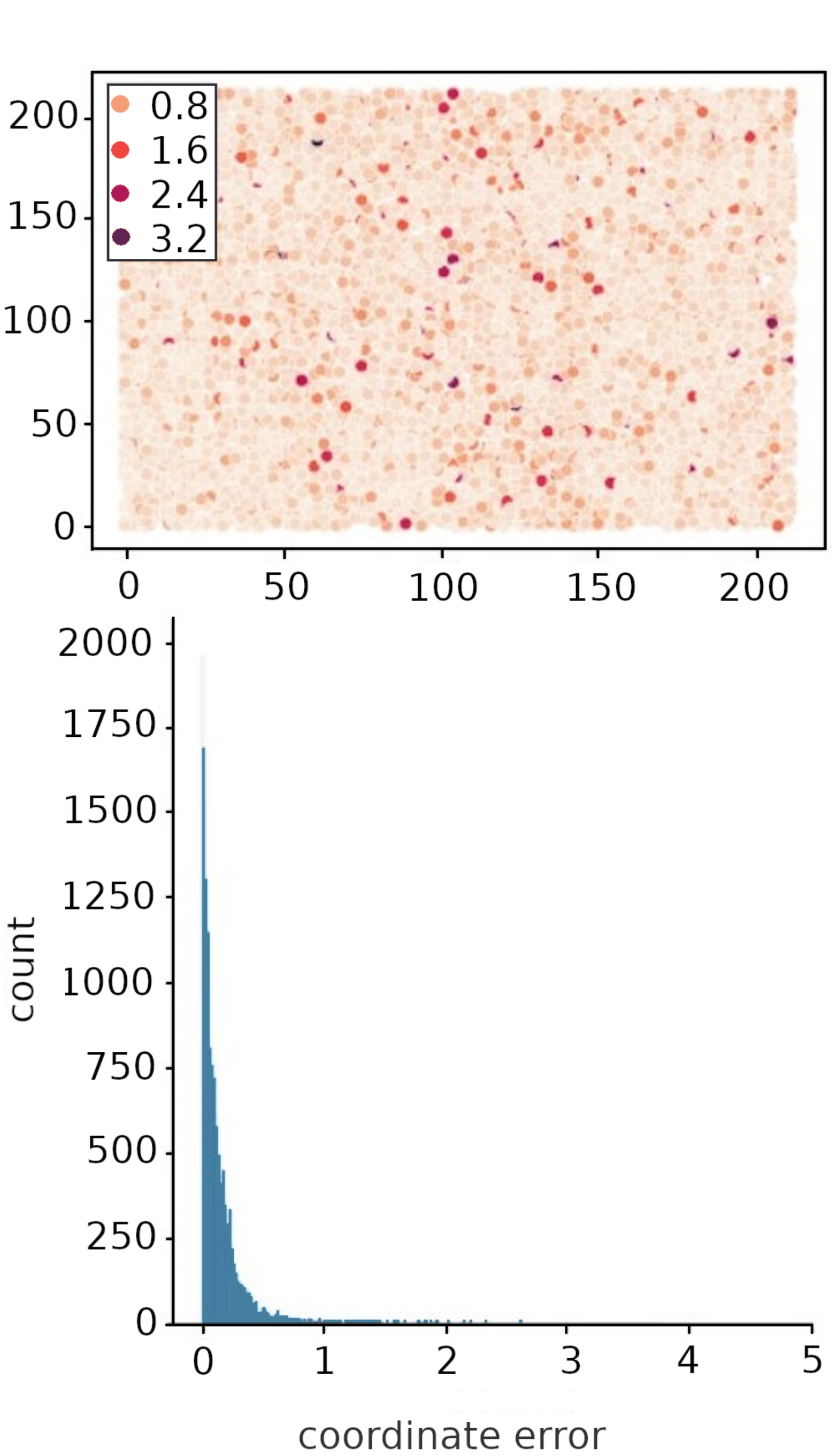}
		\caption{PILM - 8.9s}
	\end{subfigure}
	
	\caption{Coordinate errors for $N=20,000$. Scatter plot, each point is colored according to its coordinate error (row 1), and histogram of the errors (row 2).}\label{err_10}
\end{figure}

\begin{figure}[h]
	\centering
	\begin{subfigure}[b]{0.327\textwidth}
		\centering
		\includegraphics[width = \textwidth]{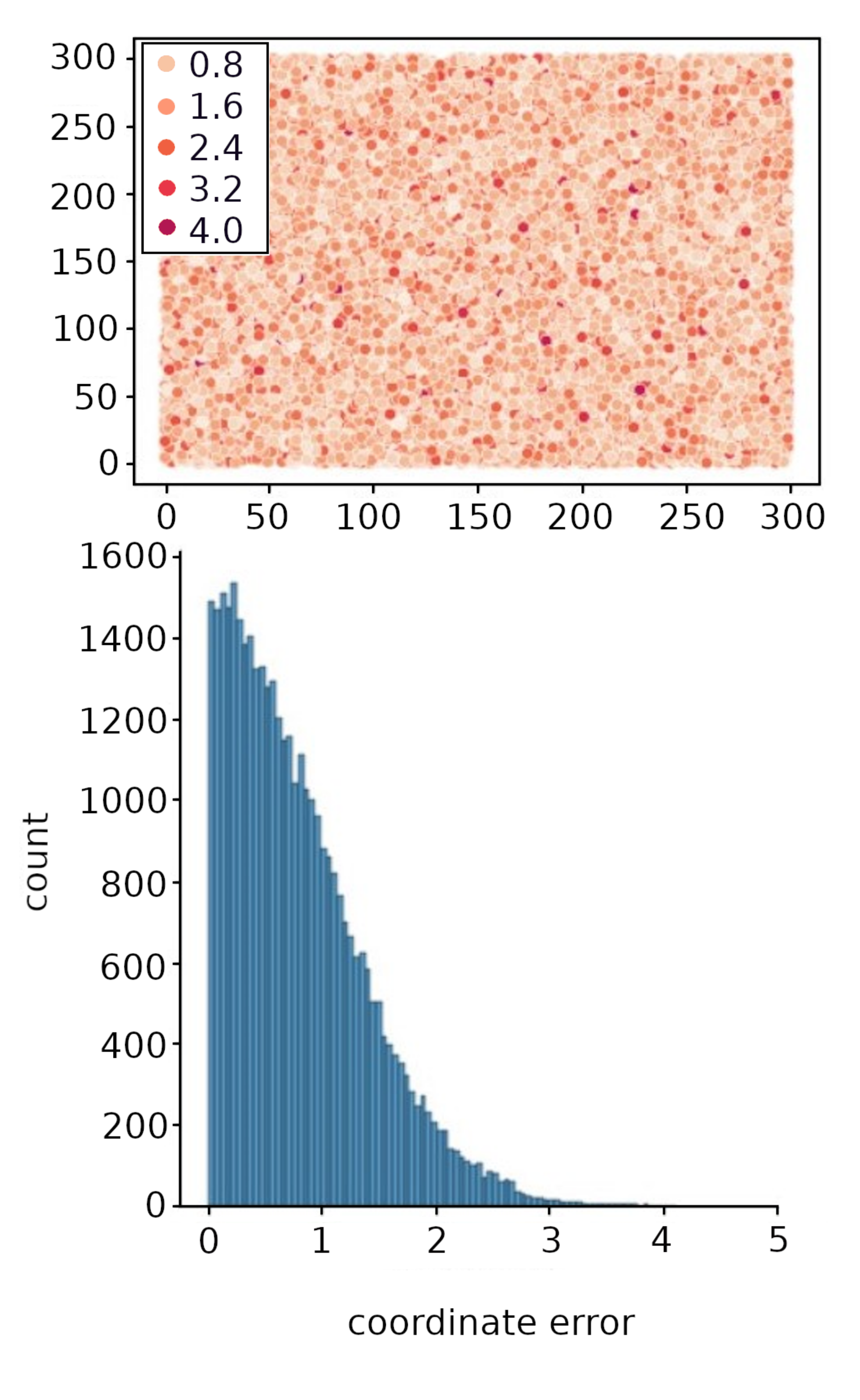}
		\caption{initial guess}
	\end{subfigure}
	\begin{subfigure}{0.327\textwidth}
		\centering
		\includegraphics[width = \textwidth]{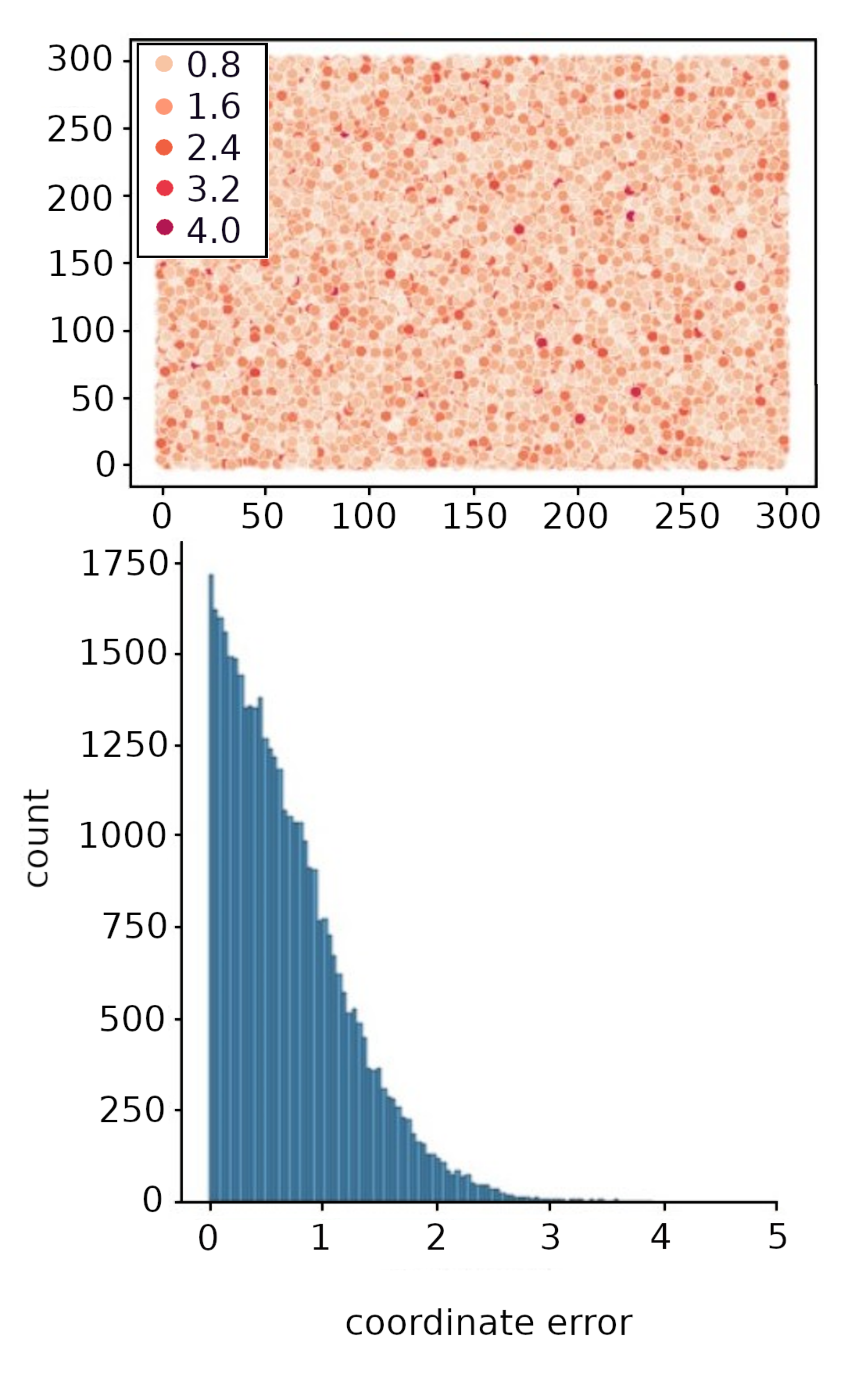}
		\caption{classical LM - 24s}
	\end{subfigure}
	\begin{subfigure}{0.327\textwidth}
		\centering
		\includegraphics[width = \textwidth]{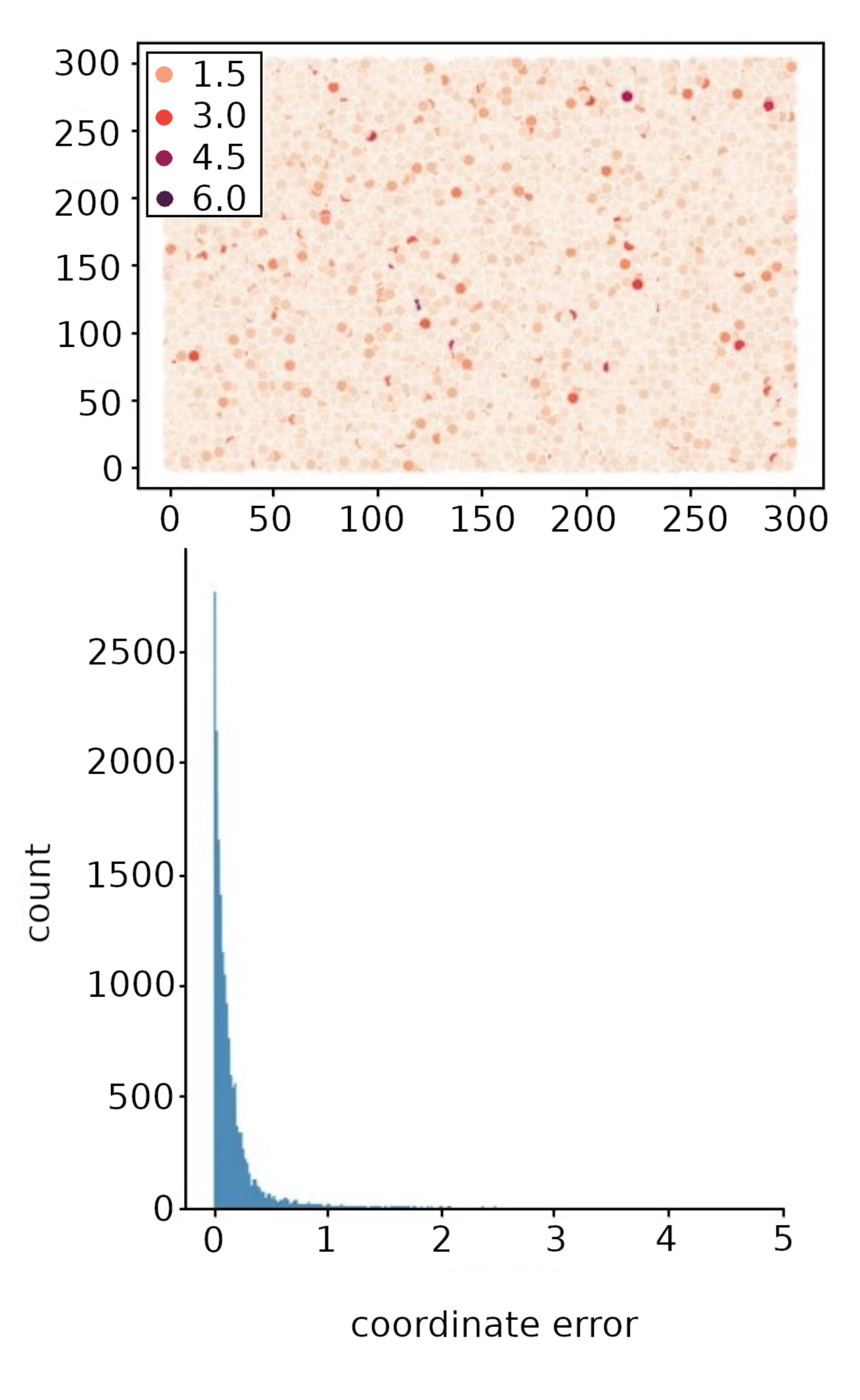}
		\caption{PILM - 15s}
	\end{subfigure}
	
	\caption{Coordinate errors for $N=40,000$. Scatter plot, each point is colored according to its coordinate error (row 1), and histogram of the errors (row 2).}\label{err_20}
\end{figure}

For the problem with ${ n=10^6}$, in { F}igure \ref{err_500} we plot the distribution plot of the coordinate error at the initial guess and at the approximate solution reached by PILM with $K=72.$ Given a vector $\xbf\in\R^N$, for $i=1,\dots,N$ the coordinate error is computed as $|x_i-x^*_i|$, where $x_i^*$ is the ground truth. 
We now compare the proposed method and classical LM in terms of accuracy reached in a given budget of computational time. We proceed as follows. We consider the computational time it takes the PILM method to reach termination, then we run the classical LM method, stopping the execution after the same amount of computational time, allowing the method to finish the current iteration. We then compare the coordinate errors at the final point for the two methods.
In Figures \ref{err_10} and \ref{err_20} we plot the results for the test cases with { 20,000} and { 40,000 } variables respectively. For each of the problems we consider the coordinate error at the initial guess (subfigure (a)), at the coordinates computed by PILM with optimal $K$ (subfigure (c)), and at the coordinates computed by the LM method in a comparable amount of computational time (subfigure (b)). \REV{For each considered case we show both a histogram of the coordinate errors and a scatter plot where the dot representing each point $(x_{2i-1},x_{2i})$ is colored depending on its relative error.} We report that for $N\geq60,000$, the LM method cannot complete the first iteration within the given budget of time.

\begin{figure}[h]
	\centering
	\includegraphics[width=0.8\linewidth]{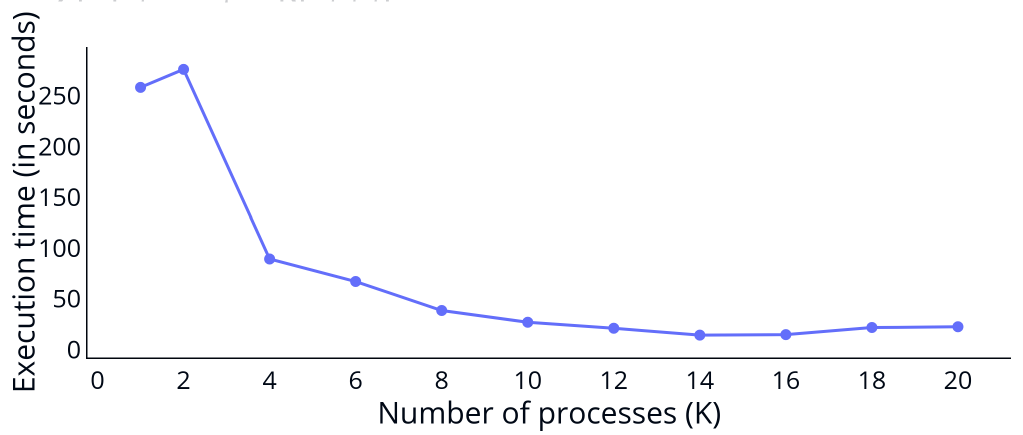}
	\caption{Execution time for different number of processors} 
	\label{fig:unbalanced}
\end{figure}

In the numerical experiments that we presented so far the partitioning of the method was well balanced, in the sense that all the subsets had approximately the same cardinality. In the following we consider the same problem with $N=40,000$ variables as in Figure 7, but the partitioning is done in such a way that some of the subsets are significantly larger than the others. We run the method for $K=1$ and $K\in\{2,4,\dots,20\}$ and consider the time necessary to reach termination. The initial guess, the termination condition and all the parameters are the same as in the previous tests. In Figure \ref{fig:unbalanced} we plot the obtained results. For completeness, in Table \ref{tab:unbalanced} we report, for each considered value of $K$, the ratio between the largest and the smallest subset in the partition. Notice that for the results presented in the previous tests, the ratio was always between 1.01 and 1.22.
\begin{table}[h]
	\centering
	\begin{tabular}{c|c|c|c|c|c|c|c|c|c|c}
		$\mathbf{K}$ & 2 & 4 & 6 &8&10&12&14&16&18&20 \\
		\hline
		\textbf{ratio} & 2.8& 3.0 & 2.5& 3.6& 3.6& 3.6& 3.6& 3.6& 3.6&3.5\\
	\end{tabular}
	\caption{Ratio between the cardinalities of the largest and smallest subset in the partition, for each value of of $K.$}
	\label{tab:unbalanced}
\end{table}

As we can see from Figure \ref{fig:unbalanced} even when the variables are partitioned in subsets of significantly different sizes the proposed method appears to be more efficient that the classical LM method, and the behavior is overall analogous to that that we observed in Figure \ref{fig:PILM1}. The only value of $K$ for which PILM seems to perform worse than classical LM method is $K=2$.\\

\begin{figure}[h]
	\centering
	\begin{subfigure}[b]{0.327\textwidth}
		\centering
		\includegraphics[width = \textwidth]{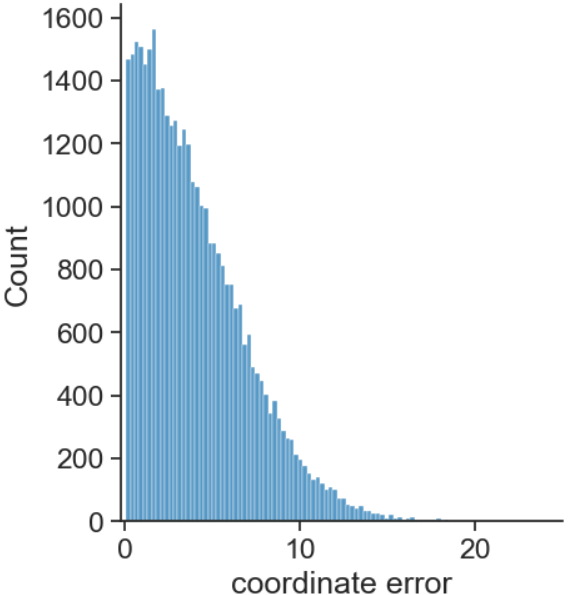}
		\caption{initial guess}
	\end{subfigure}
	\begin{subfigure}{0.327\textwidth}
		\centering
		\includegraphics[width = \textwidth]{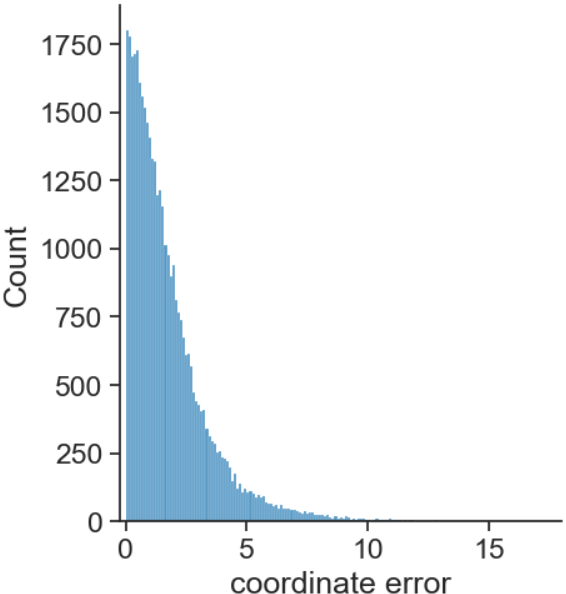}
		\caption{classical LM - 108s}
	\end{subfigure}
	\begin{subfigure}{0.327\textwidth}
		\centering
		\includegraphics[width = \textwidth]{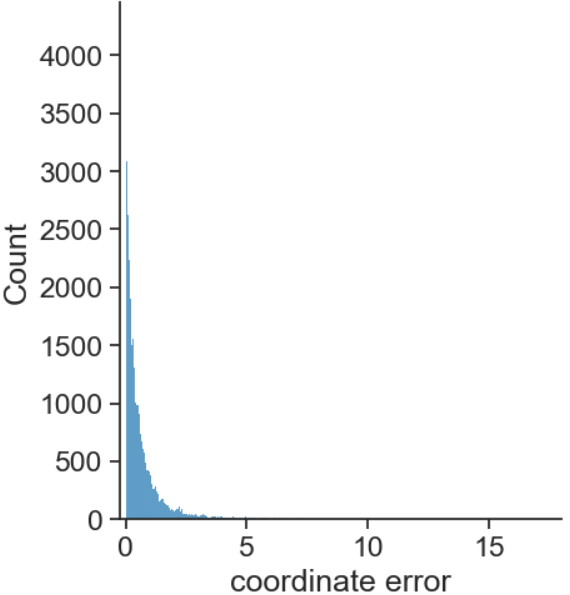}
		\caption{PILM - 100s}
	\end{subfigure}
	
	\caption{Histogram of the coordinate errors for $N=40,000$.}\label{fig:stdev5}
\end{figure}

To conclude this section, we present the results of a numerical test aimed at understanding how the quality of the initial guess influences the comparison between PILM and classical LM. We consider the same problem with $N=40,000$ variables as in Figure 7, however, the coordinate observation are generated with standard deviation 5.0 instead of 1.0. We run both PILM method with $K=14$ and classical LM method for 100 seconds. And compute the coordinate error ate the point $\xbf^k$ reached by each of the methods at the end of the last iteration. In Figure \ref{fig:stdev5} we plot the histograms of the coordinate errors. To improve the quality of the solution, for PILM method we increased the number of inner iterations to 50. All other parameters are the same as in the previous tests. From Figure \ref{fig:stdev5} we can see that also in this case PILM method achieves a significantly more accurate point, compared to the classical LM method, in the given budget of computational time.

\section{Conclusion}
We considered a generic nonlinear least squares (NLS) problem of minimizing the sum of squared per-measurement losses, where each measurement involves a few unknown coordinates. More precisely, our main interest is in large scale nearly separable NLS problems, motivated by localization-type applications such as cadastral maps refinements. The classical Levenberg Marquardt (LM) method is a widely adopted tool for NLS; however, its direct application in huge scale problems is either too costly or infeasible. In this paper, we develop an efficient parallel method based on LM dubbed PILM that harnesses the nearly separable structure of the underlying NLS problem 
for efficient parallelization. A detailed convergence analysis is provided for the proposed method. First, we prove that PILM, combined with a nonmonotone line search strategy, achieves global convergence to a stationary point of the NLS problem. For the full step size, PILM exhibits local convergence, with the convergence order depending on the choice of the parameters of the method. The achieved results for PILM hold under a standard set of assumptions, akin to those of the classical LM's theory.  
An efficient implementation of PILM is provided in a master-worker parallel compute environment, and its efficiency is demonstrated on huge scale cadastral map refinement problems.

\section*{Data Availability Statement}
The testcases used to obtain the numerical results presented in Section 4 are generated by the authors and publicly available at the following address https://cloud.pmf.uns.ac.rs/s/GaSNnns9fdJeXqD. The code is available at https://github.com/lidijaf/PILM.

\bibliographystyle{plain}
\bibliography{references}

\end{document}